\documentclass[10pt]{amsart}
\usepackage{amssymb}
\usepackage{amscd}

\theoremstyle{definition}
\newtheorem{thm}{Theorem}[section]
\newtheorem{lem}[thm]{Lemma}
\newtheorem{prp}[thm]{Proposition}
\newtheorem{dfn}[thm]{Definition}
\newtheorem{cor}[thm]{Corollary}

\newtheorem{rmk}[thm]{Remark}

\newtheorem{exa}[thm]{Example}

\newcommand{\beq}{\begin{equation}}
\newcommand{\eeq}{\end{equation}}
\newcommand{\beqr}{\begin{eqnarray*}}
\newcommand{\eeqr}{\end{eqnarray*}}
\newcommand{\bal}{\begin{align*}}
\newcommand{\eal}{\end{align*}}
\newcommand{\bei}{\begin{itemize}}
\newcommand{\eei}{\end{itemize}}

\newcommand{\af}{\alpha}
\newcommand{\bt}{\beta}
\newcommand{\gm}{\gamma}
\newcommand{\dt}{\delta}
\newcommand{\ep}{\varepsilon}
\newcommand{\zt}{\zeta}

\newcommand{\ch}{\chi}
\newcommand{\io}{\iota}
\newcommand{\te}{\theta}
\newcommand{\ld}{\lambda}

\newcommand{\ph}{\varphi}
\newcommand{\ps}{\psi}
\newcommand{\rh}{\rho}
\newcommand{\om}{\omega}
\newcommand{\ta}{\tau}

\newcommand{\Gm}{\Gamma}

\newcommand{\Q}{{\mathbf{Q}}}
\newcommand{\Z}{{\mathbf{Z}}}
\newcommand{\R}{{\mathbf{R}}}
\newcommand{\C}{{\mathbf{C}}}
\newcommand{\N}{{\mathbf{N}}}

\newcommand{\id}{{\mathrm{id}}}

\newcommand{\dist}{{\mathrm{dist}}}

\newcommand{\Ad}{{\mathrm{Ad}}}

\newcommand{\Aut}{{\mathrm{Aut}}}

\newcommand{\Ker}{{\mathrm{Ker}}}

\newcommand{\Mi}{M_{\infty}}

\newcommand{\andeqn}{\,\,\,\,\,\, {\mbox{and}} \,\,\,\,\,\,}

\newcommand{\ts}[1]{{\textstyle{#1}}}

\newcommand{\ca}{C*-algebra}
\newcommand{\ct}{continuous}
\newcommand{\pj}{projection}

\newcommand{\nbhd}{neighborhood}

\newcommand{\hm}{homomorphism}
\newcommand{\fd}{finite dimensional}
\newcommand{\wolog}{without loss of generality}
\newcommand{\Wolog}{Without loss of generality}

\newcommand{\ifo}{if and only if}

\newcommand{\mops}{mutually orthogonal \pj s}

\newcommand{\hme}{homeomorphism}

\newcommand{\tgca}{transformation group \ca}

\newcommand{\cp}{crossed product}

\newcommand{\cfn}{continuous function}
\newcommand{\hsa}{hereditary subalgebra}

\newcommand{\mvnt}{Murray-von Neumann equivalent}

\newcommand{\tRp}{tracial Rokhlin property}

\newcommand{\sfsuca}{stably finite simple unital \ca}

\title[Noncommutative Furstenberg transformations]{Furstenberg
    transformations on irrational rotation algebras}

\author{Hiroyuki Osaka}
\author{N.\  Christopher Phillips}

\date{2 September 2004}

\address{Department of Mathematics, Ritsumeikan University,
       Kusatsu, Shiga, 525-8577, Japan.}
\email[]{osaka@se.ritsumei.ac.jp}

\address{Department of Mathematics, University  of Oregon,
       Eugene OR 97403-1222, USA.}
\email[]{ncp@darkwing.uoregon.edu}

\thanks{
Research partially supported by
JSPS Grant for Scientific Research No.\  14540217(c)(1),
and by NSF grants DMS-0070776 and DMS-0302401.}

\begin{document}

\begin{abstract}
We introduce a general class of automorphisms of rotation algebras,
the noncommutative Furstenberg transformations.
We prove that
fully irrational noncommutative Furstenberg transformations
have the tracial Rokhlin property,
which is a strong form of outerness.
We conclude that crossed products by these automorphisms
have stable rank one, real rank zero,
and order on projections determined by traces
(Blackadar's Second Fundamental Comparability Question).

We also prove that several classes of simple quotients of the
C*-algebras of discrete subgroups of five dimensional
nilpotent Lie groups, considered by Milnes and Walters,
are crossed products of simple C*-algebras
(C*-algebras of minimal ordinary Furstenberg transformations)
by automorphisms which have the tracial Rokhlin property.
It follows that these algebras also have
stable rank one, real rank zero,
and order on projections determined by traces.
\end{abstract}

\maketitle

\setcounter{section}{-1}

\section{Introduction}\label{Sec:Intro}

\indent
Furstenberg introduced in~\cite{Fr}
a family of \hme s of $S^1 \times S^1,$
now called Furstenberg transformations.
They have the form
\[
h_{\gm, d, f} (\zt_1, \zt_2)
 = (e^{2 \pi i \gm} \zt_1, \, \exp (2 \pi i f (\zt_1)) \zt_1^d \zt_2),
\]
with fixed $\gm \in \R,$ $d \in \Z,$ and $f \colon S^1 \to \R$ \ct.
For $\gm \not\in \Q$ and $d \neq 0,$
Furstenberg proved that $h_{\gm, d, f}$ is minimal.
These \hme s, and higher dimensional analogs
(which also appear in~\cite{Fr}),
have attracted significant interest in operator algebras
(see, for example, \cite{Pc86}, \cite{Ji}, \cite{Kd}, and~\cite{RM})
and in dynamics (see, for example, \cite{ILR} and~\cite{Rh}).

For any $\te \in \R,$
the formula for the automorphism $f \mapsto f \circ h_{\gm, d, f}$
of $C (S^1 \times S^1)$ also defines an automorphism of the
rotation algebra $A_{\te}.$
Taking the generators of $A_{\te}$
to be unitaries $u$ and $v$ satisfying $v u = e^{2 \pi i \te} u v,$
we obtain an automorphism $\af_{\te, \gm, d, f}$ of $A_{\te}$ such that
\[
\af_{\te, \gm, d, f} (u) = e^{2 \pi i \gm} u
\andeqn \af_{\te, \gm, d, f} (v) = \exp (2 \pi i f (u)) u^d v.
\]
We call it a Furstenberg transformation on $A_{\te}.$
When $\te \not\in \Q,$
it is the most general automorphism $\af$ of $A_{\te}$
for which $\af (u)$ is a scalar multiple of $u$
(Proposition~\ref{MostGen}).
Several special cases of these have been considered previously
(see the discussion after Proposition~\ref{MostGen}),
but, as far as we know, not the general case.

In this paper,
we prove that when $1, \te, \gm$ are linearly independent over $\Q,$
then $\af_{\te, \gm, d, f}$ has the \tRp\  introduced in~\cite{OP1}.
In particular, all nontrivial powers of $\af_{\te, \gm, d, f}$
are outer in a strong sense.
Combining this result with that of~\cite{OP1},
we see that $C^* (\Z, \, A_{\te}, \, \af_{\te, \gm, d, f})$
has stable rank one and real rank zero,
and that the order on \pj s over this algebra is determined by traces
(Blackadar's Second Fundamental Comparability Question).

Other cases, in which $1, \te, \gm$ are linearly dependent over $\Q,$
or even in which $\te,$ $\gm,$ or both are rational,
are also interesting, but the analysis becomes more difficult.
Sometimes the automorphism has the \tRp\  %
(see Theorem~\ref{MainLemma2}),
sometimes the automorphism does not have the \tRp\  but the
crossed product is simple anyway (see Example~\ref{NonUErg}),
and sometimes the crossed product is not simple
(such as when $\te$ and $\gm$ are both rational).
It seems to be an interesting question to understand the structure
of the crossed products in these cases,
but we leave it for future work.

The method of proof applies to other examples as well.
In a series of
papers~\cite{MW1}, \cite{MW2}, \cite{MW3}, \cite{MW4},~\cite{Wl},
Milnes and Walters have studied the simple quotients of the
\ca s of certain discrete subgroups of nilpotent Lie groups
of dimension up to five.
These algebras are a kind of generalization of the irrational
rotation algebras, which occur when the Lie group is the three
dimensional Heisenberg group.
Among the algebras in these papers, the least is known about
the algebras $A_{\te}^{5, 3}$ and $A_{\te}^{5, 6}$ of~\cite{MW2}
and the algebras $A_{\te}^{5, 3} (\af, \bt, \gm, \dt, \ep)$
of~\cite{MW4}.
We show that each of these is the crossed product of a simple \ca\  %
(the \ca\  of an ordinary minimal Furstenberg transformation
on $S^1 \times S^1$)
by an automorphism with the \tRp.
Again, it follows that these algebras
have stable rank one and real rank zero,
and that the order on \pj s over them is determined by traces.

This paper is organized as follows.
In Section~\ref{Sec:Prelim},
we define Furstenberg transformations on rotation algebras,
and prove some basic properties.
Among other things, we compute the ordered K-theory of the crossed
product when $1, \te, \gm$ are linearly independent over $\Q,$
and we identify several cases
in which one can prove by elementary means
that crossed products of two different rotation algebras by
two different Furstenberg transformations are isomorphic.

In Section~\ref{Sec:Outer}, we derive an easy to verify
criterion for the \tRp.
Let $A$ be a simple unital \ca\  with tracial rank zero in the sense
of \cite{LnTAF},~\cite{LnTTR},
and with a unique tracial state.
Then an automorphism of $A$ has the \tRp\  \ifo\  the corresponding
automorphism of the type~II$_1$ factor obtained from the tracial
state via the Gelfand-Naimark-Segal construction is aperiodic.
This is a kind of generalization of the equivalence of
Conditions~(i) and~(iii) in Theorem~2.1 of~\cite{Ks3},
and we essentially follow the proof there.
In Section~\ref{Sec:OutPf}, we apply this result to
$\af_{\te, \gm, d, f}$ in the case that
$1, \te, \gm$ are linearly independent over $\Q,$
and also in the case
$\te \not\in \Q,$ $\gm = 0,$ $d \neq 0,$ and $f$ constant.
Section~\ref{Sec:MW} contains the results on simple quotients of
discrete subgroups of five dimensional nilpotent Lie groups.

Both authors would like to thank Masaki Izumi
for valuable discussions,
and the first author would also like to thank Kazunori Kodaka.

\section{Furstenberg transformations on rotation
   algebras}\label{Sec:Prelim}

\indent
For $\te \in \R$ let $A_{\te}$ be the usual rotation algebra,
generated by unitaries $u$ and $v$ satisfying
$v u = e^{2 \pi i \te} u v.$

\begin{dfn}\label{FDfn}
Let $\te, \, \gm \in \R,$ let $d \in \Z,$
and let $f \colon S^1 \to \R$ be a \cfn.
The {\emph{Furstenberg transformation}} on $A_{\te}$
determined by $(\te, \gm, d, f)$ is the automorphism
$\af_{\te, \gm, d, f}$ of $A_{\te}$ such that
\[
\af_{\te, \gm, d, f} (u) = e^{2 \pi i \gm} u
\andeqn \af_{\te, \gm, d, f} (v) = \exp (2 \pi i f (u)) u^d v.
\]
\end{dfn}

We check below that $\af_{\te, \gm, d, f}$ is an automorphism.
It might be called a noncommutative Furstenberg transformation.
Note that the parameter $\te$ does not appear in the formulas;
its only role is to specify the algebra on which the automorphism
acts.

Special cases of these have been studied previously.
See the discussion after Proposition~\ref{MostGen}.

\begin{lem}\label{ForDfn}
Let $\te, \, \gm \in \R,$ let $d \in \Z,$
and let $f \colon S^1 \to \R$ be a \cfn.
Then there is a unique automorphism given by the formulas in
Definition~\ref{FDfn}.
\end{lem}

\begin{proof}
One checks easily that
\[
u_0 = e^{2 \pi i \gm} u \andeqn
v_0 = \exp (2 \pi i f (u)) u^d v
\]
are unitaries which satisfy the commutation relation
$v_0 u_0 = e^{2 \pi i \te} u_0 v_0.$
Therefore there is a \hm\  $\af_{\te, \gm, d, f}$
which sends $u$ to $u_0$ and $v$ to $v_0.$
Uniqueness is obvious.
To prove that it is an automorphism, we exhibit an inverse.
With
$g (\zt) = - f \left( e^{- 2 \pi i \gm} \zt \right) + d \gm$
for $\zt \in S^1,$
one can check that $\af_{\te, \, - \gm, \, - d, \, g}$
is a two sided inverse for $\af_{\te, \gm, d, f}.$
\end{proof}

We are interested in the crossed products
$C^* (\Z, \, A_{\te}, \, \af_{\te, \gm, d, f})$ by these automorphisms.
There is some immediate redundancy,
which we exhibit in the next three lemmas.

\begin{lem}\label{ChangeOfVar}
Let $\te, \, \gm \in \R,$ let $d \in \Z,$
and let $f \colon S^1 \to \R$ be a \cfn.
Then
\[
C^* (\Z, \, A_{\te}, \, \af_{\te, \gm, d, f})
  \cong C^* (\Z, \, A_{\gm}, \, \af_{\gm, \, \te, \, - d, \, - f}).
\]
\end{lem}

\begin{proof}
Let $w$ be the canonical implementing unitary of the crossed product
\[
C^* (\Z, \, A_{\te}, \, \af_{\te, \gm, d, f}).
\]
Then $C^* (\Z, \, A_{\te}, \, \af_{\te, \gm, d, f})$ is the universal
\ca\  generated by three unitaries $u, v, w,$ subject to the relations
\[
v u = e^{2 \pi i \te} u v, \,\,\,\,\,\,
w u w^* = e^{2 \pi i \gm} u, \andeqn
w v w^* = \exp (2 \pi i f (u)) u^d v.
\]
We take adjoints in the third relation, exchange the first two,
and rearrange algebraically, to get
\[
w u = e^{2 \pi i \gm} u w, \,\,\,\,\,\,
v u v^* = e^{2 \pi i \te} u, \andeqn
v w v^* = \exp (- 2 \pi i f (u)) u^{-d} w.
\]
The universal \ca\  generated by unitaries satisfying these
relations is
\[
C^* (\Z, \, A_{\gm}, \, \af_{\gm, \, \te, \, - d, \, - f}),
\]
with $v$ playing the role of the
implementing unitary of the crossed product.
\end{proof}

\begin{lem}\label{Conj}
Let $\te, \, \gm \in \R,$ let $d \in \Z,$
and let $f \colon S^1 \to \R$ be a \cfn.
Let $k, l \in \Z.$
Set $g (\zt) = f \left( e^{2 \pi i l \te} \zt \right) + (l d - k) \te.$
Then
\[
\Ad (u^k v^l) \circ \af_{\te, \gm, d, f}
  = \af_{\te, \, \gm + l \te, \, d, \, g}.
\]
Moreover, the crossed products by $\af_{\te, \gm, d, f}$
and $\af_{\te, \, \gm + l \te, \, d, \, g}$ are isomorphic.
\end{lem}

\begin{proof}
The first statement is a computation.
The second is then a special case of the isomorphism
$C^* (\Z, \, A, \, \Ad(u) \circ \af) \cong C^* (\Z, A, \af)$
for any unitary $u.$
Lacking a suitable reference for this isomorphism,
we outline the proof here.
The actions of $\Z$ generated by $\af$ and $\Ad(u) \circ \af$
are exterior equivalent, in the sense of 8.11.3 of~\cite{Pd1},
via the cocycle defined by $u_n = u \af (u) \cdots \af^{n - 1} (u)$
for $n \geq 0$ and $u_{- n} = \af^{- n} (u_n^*)$ for $n > 0.$
The crossed products are then isomorphic as in the proof of
Theorem 2.8.3(5) of~\cite{Ph1}.
\end{proof}

\begin{lem}\label{Conj2}
Let $\te, \, \gm \in \R,$ let $d \in \Z \setminus \{ 0 \},$
and let $f \colon S^1 \to \R$ be a \cfn.
Let $r \in \R,$
and set $g (\zt) = f \left( e^{2 \pi i r / d} \zt \right) + r.$
Then the automorphism $\bt$ of $A_{\te},$ determined by
$\bt (u) = \exp (2 \pi i r / d) u$ and $\bt (v) = v,$
satisfies
\[
\bt \circ \af_{\te, \gm, d, f} \circ \bt^{-1}
  = \af_{\te, \, \gm, \, d, \, g}.
\]
Moreover, the crossed products by $\af_{\te, \gm, d, f}$
and $\af_{\te, \, \gm, \, d, \, g}$ are isomorphic.
\end{lem}

\begin{proof}
The first statement is a computation.
The second statement is immediate from the first.
\end{proof}

In this lemma, one does not get anything new by
taking $\bt (v) = \exp (2 \pi i s) v.$

As we now show, when $\te \in \R \setminus \Q,$
automorphisms of the form $\af_{\te, \gm, d, f}$
are the most general automorphisms of $A_{\te}$ which send
$u$ to a scalar multiple of $u.$

\begin{prp}\label{MostGen}
Let $\te \in \R \setminus \Q$ and let $\gm \in \R.$
Let $\af \in \Aut (A_{\te})$ be an automorphism
such that $\af (u) = e^{2 \pi i \gm} u.$
Then there exist $d \in \Z$ and a \cfn\  $f \colon S^1 \to \R$
such that $\af = \af_{\te, \gm, d, f}.$
\end{prp}

\begin{proof}
We have
\[
u \af (v) u^* = [e^{2 \pi i \gm} u] \af (v) [e^{2 \pi i \gm} u]^*
  = \af (u v u^*) = e^{- 2 \pi i \te} \af (v).
\]
Therefore $u$ commutes with $\af (v) v^*.$
Since $\te$ is irrational, $C^* (u)$ is a maximal commutative
subalgebra of $A_{\te}$ by Corollary 3.3.3 on page~79 of~\cite{Tm}.
(See the previous page for the notation.)
So there is a \cfn\  $g \colon S^1 \to S^1$ such that
$\af (v) v^* = g (u).$
We may write $g (\zt) = \exp (2 \pi i f (\zt)) \zt^d$
for suitable $d \in \Z$ (the degree of $g$)
and a suitable \cfn\  $f \colon S^1 \to \R.$
Then $\af (v) = \exp (2 \pi i f (u)) u^d v,$
whence $\af = \af_{\te, \gm, d, f}.$
\end{proof}

When $\te = 0,$ $d \neq 0,$ and $\gm \in \R \setminus \Q,$
then $\af_{\te, \gm, d, f}$ is a minimal
Furstenberg transformation on the torus,
as first introduced in Section~2 of~\cite{Fr}.
Crossed products by these have been well studied,
sometimes only in special cases such as $f$ constant or $d = 1,$
and sometimes in considerably greater generality,
such as in higher dimensions.
See, for example, Section~4 of~\cite{Pc86}, \cite{Rh}, \cite{Kd},
Example~4.9 of~\cite{Ph7}, and~\cite{RM}.
When $f$ is smooth, this \hme\  is uniquely ergodic by
Theorem~2.1 of~\cite{Fr}, and
it follows from~\cite{LP2}
(also see the survey article~\cite{LP1})
that the \cp\  is classifiable.
More recently,
classifiability has been proved whenever
$f$ is uniquely ergodic~\cite{LhP}.
In particular, the \cp\  does not depend on the choice
of $f,$ and when $d = 1$ it is an AT~algebra.
For arbitrary \ct\  $f,$ Theorem~2 (in Section~4) of~\cite{ILR}
shows that the Furstenberg transformation
need not be uniquely ergodic.
(See the discussion in Example~\ref{NonUErg}.)

The case $d = 1$ and $f$ equal to the constant function~$0$
has been studied by Packer \cite{Pc87},~\cite{Pc88}.
She classified the crossed products up to isomorphism
and Morita equivalence, and computed the ordered K-theory.
In the case that at least one of $\te$ and $\gm$ is not rational,
she also proved that the \pj s satisfy cancellation.

The case $\gm = 0,$ $d = 1,$ and $\te \in \R \setminus \Q$ has been
considered by Kodaka and Osaka~\cite{KO}.
They have shown that $\af_{\te, 0, 1, f}$ is outer.
Lemma~\ref{ChangeOfVar} shows that in this case the \cp\  is
isomorphic to the \tgca\  of a minimal
Furstenberg transformation on the the ordinary torus.
Using Lemmas~\ref{Conj} and~\ref{ChangeOfVar}, one can
show that if $\te \in \R \setminus \Q$ but
$1, \te, \gm$ are linearly dependent over $\Q,$
then
\[
C^* (\Z, \, A_{\te}, \, \af_{\te, \gm, d, f})
  \cong C^* (\Z, \, A_{\te_0}, \, \af_{\te_0, \gm_0, d_0, f_0})
\]
with $\te_0 \in \Q,$ that is, the \cp\  is isomorphic to the
\cp\  of a Furstenberg transformation on a rational rotation algebra.
Although we will not consider such \cp s in this paper,
it is plausible that one can relate them to the \tgca s of
ordinary Furstenberg transformations using Morita equivalence.

We compute the K-theory of the crossed product by a
Furstenberg transformation on $A_{\te}$
when $\te$ and $\gm$ are not both rational.
The computation of the groups is valid without this restriction.
The statement about the range of a tracial state should be as well,
but we do not investigate the rational case further here.

\begin{lem}\label{KThyComp}
Let $\te, \, \gm \in \R,$ let $d \in \Z,$
and let $f \colon S^1 \to \R$ be a \cfn.
Let $B = C^* (\Z, \, A_{\te}, \, \af_{\te, \gm, d, f}).$
Suppose $\te$ and $\gm$ are not both rational and $d \neq 0.$
Then
\[
K_0 ( B ) \cong \Z^3
\andeqn
K_1 ( B ) \cong \Z^3 \oplus \Z / d \Z.
\]
Moreover, if $\ta$ is any tracial state on $B,$ then
$\ta_* (K_0 (B)) = \Z + \te \Z + \gm \Z.$
\end{lem}

\begin{proof}
By Lemma~\ref{ChangeOfVar}, \wolog\  $\te \not\in \Q.$
Let $\af = \af_{\te, \gm, d, f}.$
\Wolog\  $0 < \te < 1.$
Let $\io \colon A_{\te} \to B$ be the inclusion map.
The composition $\ta \circ \io$ is the unique tracial state on $A_{\te},$
and we denote it again by $\ta.$

We use the Pimsner-Voiculescu exact sequence
\cite{PV}, which here takes the form
\[
\begin{CD}
K_0 (A_{\te}) @>{\id - \af^{-1}_*}>> K_0 (A_{\te}) @>{\io_*}>> K_0 (B)\\
@A{\mathrm{exp}}AA & &  @VV{\partial}V\\
K_1 (B)  @<<{\io_*}< K_1 (A_{\te}) @<<{\id - \af^{-1}_*}< K_1 (A_{\te}).
\end{CD}
\]
Recall that $K_0 (A_{\te}) \cong \Z^2$ with generators $[1]$
and $[p],$ where $p \in A_{\te}$ is a \pj\  with $\ta (p) = \te.$
Also, recall that $K_1 (A_{\te}) \cong \Z^2$ with generators $[u]$
and $[v],$ the classes of the standard unitary generators.
The induced map $\af_*$ on $K_0 (A_{\te})$ is the identity,
since $K_0 (A_{\te})$ has no nontrivial automorphisms
preserving order and $[1]$ when $\te \not\in \Q.$
We have $\left[ e^{2 \pi i \gm} u \right] = [u]$ and
$\left[ \exp (2 \pi i f (u)) u^d v \right] = d [u] + [v],$
so $\af_*$ on $K_1 (A_{\te})$ has the matrix
{\scriptsize{$
 \left( \begin{array}{cc} 1 & d \\ 0 & 1 \end{array} \right) $}}.
Therefore the upper left horizontal map in the sequence above is
zero and the lower right horizontal map is
{\scriptsize{$
 \left( \begin{array}{cc} 0 & d \\ 0 & 0 \end{array} \right) $}}.
The computation of the groups follows easily.

It remains to determine the range of $\ta_*.$
If $\gm = 0$ then the result on the range of the tracial state follows from
Lemma~\ref{ChangeOfVar} and Example~4.9 of~\cite{Ph7}.
If $\gm$ is rational, we may use Lemma~\ref{Conj} to replace
$\gm$ by $\gm + \te.$
Accordingly, we may assume $\gm \not\in \Q.$
We may further assume $0 < \gm < 1.$

Let $u_0 \in C (S^1)$ be the standard generating unitary,
and let $\ph \colon C (S^1) \to A_{\te}$ be the \hm\  determined by
$\ph (u_0) = u.$
Then $\ph$ is equivariant for the automorphisms
$u_0 \to e^{2 \pi i \gm} u_0$ of $C (S^1)$ and $\af$ of $A_{\te}.$
Therefore there is a \hm\  $\ps \colon A_{\gm} \to B$ of the \cp s.
Taking suitable portions of the Pimsner-Voiculescu exact sequences
for these \cp s, we get the following commutative diagram with exact
rows:
\[
\begin{CD}
0 @>{}>> K^0 (S^1) @>{(\io_0)_*}>> K_0 (A_{\gm})
                     @>{\partial_0}>> K^1 (S^1) @>{}>> 0 \\
& &    @V{\ph_*}VV   @V{\ps_*}VV   @V{\ph_*}VV    \\
0 @>{}>> K_0 (A_{\te}) @>{\io_*}>> K_0 (B)
     @>{\partial}>> K_1 (A_{\te}) @>{\id - \af^{-1}_*}>> K_1 (A_{\te}).
\end{CD}
\]
As for $A_{\te},$ the composition $\ta \circ \ps$ is the
unique tracial state on $A_{\gm},$ and we write simply $\ta.$
Further, $K_0 (A_{\gm}) \cong \Z^2$ with generators $[1]$
and $[q],$ where $q \in A_{\gm}$ is a \pj\  with $\ta (q) = \gm.$
Moreover, $\partial_0 ([q]) = [u_0].$

We claim that $K_0 (B)$ is generated by
$[1],$ $\io_* ([p]),$ and $\ps_* ([q]).$
Since these have trace $1,$ $\te,$ and $\gm,$ the result will follow.
Let $\om \in K_0 (B).$
Then
\[
\partial (\om) \in \Ker (\id - \af^{-1}_*) = \Z \cdot [u].
\]
Choose $m \in \Z$ such that $\partial (\om) = m [u].$
Then
\[
\partial (m  \ps_* ([q])) = m \cdot \ph_* \circ \partial_0 ([q])
  = m \ph_* ([u_0]) = m [u] = \partial (\om).
\]
Therefore $\partial (m  \ps_* ([q]) - \om) = 0,$ so there are
$k, l \in \Z$ such that
$m  \ps_* ([q]) - \om = \io_* (k [1] + l [p]).$
Thus $\om = k [1] + l \io_* ([p]) + m  \ps_* ([q]),$ as desired.
\end{proof}

To further analyze the crossed products by the
Furstenberg transformations $\af_{\te, \gm, d, f}$
with $1, \te, \gm$ are linearly independent over $\Q,$
we need to show that these automorphisms have the \tRp.
In the next section, we give a general method for proving
this property.

\section{The tracial Rokhlin property in terms of trace
           norms}\label{Sec:Outer}

\indent
We recall the definition of the \tRp\  (Definition~1.1 of~\cite{OP1}).

\begin{dfn}\label{TRPDfn}
Let $A$ be a \sfsuca\   and let $\af \in \Aut (A).$
We say that $\af$ has the {\emph{tracial Rokhlin property}}
if for every finite set $F \subset A,$ every $\ep > 0,$
every $n \in \N,$
and every nonzero positive element $x \in A,$
there are \mops\  $e_0, e_1, \ldots, e_n \in A$ such that:
\begin{itemize}
\item[(1)]
$\| \af (e_j) - e_{j + 1} \| < \ep$ for $0 \leq j \leq n - 1.$
\item[(2)]
$\| e_j a - a e_j \| < \ep$ for $0 \leq j \leq n$ and all $a \in F.$
\item[(3)]
With $e = \sum_{j = 0}^{n} e_j,$ the \pj\  $1 - e$ is \mvnt\  to a
\pj\  in the \hsa\  of $A$ generated by $x.$
\end{itemize}
\end{dfn}

The definition requires no condition on $\af (e_n).$

The main result of this section is a criterion for
the \tRp\  for an automorphism $\af$ of a simple
unital \ca\  $A$ with tracial rank zero and unique tracial state,
namely that $\af$ should induce an aperiodic automorphism of the
type II$_1$ factor obtained as the weak operator closure of
$A$ in the Gelfand-Naimark-Segal representation associated with $\ta.$
It is a kind of generalization of the equivalence of
Conditions~(i) and~(iii) in Theorem~2.1 of~\cite{Ks3}.
We follow the outline in Section~4 of~\cite{Ks3},
but we give the details omitted there.
Some of the lemmas here are of independent interest,
such as the approximate supremum of finitely many \pj s in
a \ca\  with real rank zero (Lemma~\ref{Z3}).
In Theorem~\ref{TraceVersion},
we also show, under the same hypotheses as above,
that the \tRp\  is equivalent to an analog of the Rokhlin property
in which operator norm estimates are replaced by trace norm estimates.

Since the following computation will be used several times,
we state it separately.

\begin{lem}\label{Calc}
Let $A$ be a \ca, and let $a \in A$ satisfy
$0 \leq a \leq 1.$
Then $\| a b - b \| \leq [2 \| a b b^* - b b^* \| ]^{1/2}$
for every $b \in A.$
\end{lem}

\begin{proof}
We have
\[
\| a b - b \|^2
   \leq  \| a b b^* - b b^* \| \cdot\| a^* \| + \| a b b^* - b b^* \|
   \leq 2 \| a b b^* - b b^* \|.
\]
Take square roots.
\end{proof}

The following lemma is an approximate version for
\ca s with real rank zero of the projection onto the
closed span of finitely many subspaces.
We get the correct size control on the \pj,
but it only approximately does the correct thing on most of the
subspaces.

\begin{lem}\label{Z3}
Let $A$ be a \ca\  with real rank zero.
Let $q_0, q_1, \ldots, q_n \in A$ be \pj s, and let $\ep > 0.$
Then there exists a \pj\  $e \in A$ such that:
\begin{itemize}
\item[(1)]
$q_0 \leq e.$
\item[(2)]
$\| e q_k - q_k \| < \ep$ for $1 \leq k \leq n.$
\item[(3)]
$e \precsim q_0 \oplus q_1 \oplus \cdots \oplus q_n$ in $M_{n + 1} (A).$
\end{itemize}
\end{lem}

\begin{proof}
We first prove this when $n = 1.$
Set $p = q_0$ and $q = q_1.$
Since $A$ has real rank zero,
there is a \pj\  $r \in {\overline{ (1 - p) q A q (1 - p) }}$
such that
\[
\| r (1 - p) q (1 - p) - (1 - p) q (1 - p) \| < \ts{\frac{1}{2}} \ep^2.
\]
Then
$\| r (1 - p) q - (1 - p) q \| < \ep$ by Lemma~\ref{Calc}.
Set $e = p + r.$
Clearly $e \geq p,$ which is~(1).
We have
\begin{align*}
\| e q - q \|
 & = \| e p q + e (1 - p) q  - q \|
   = \| p q + r (1 - p) q  - q \|   \\
 & = \| r (1 - p) q - (1 - p) q \| < \ep,
\end{align*}
which is~(2).
Finally, it follows from Lemma~4.1 of~\cite{OP1} that
$r$ is \mvnt\  to a \pj\  in
${\overline{ q (1 - p) A (1 - p) q }} \subset q A q.$
So $e \precsim p \oplus q,$ which is~(3).

The general case is proved by induction.
If the statement holds for $n,$
and $q_0, q_1, \ldots, q_{n + 1} \in A$ are \pj s,
then we find $e_0$ for $q_0, q_1, \ldots, q_n \in A$ by
the induction hypothesis and $e$ for the \pj s $e_0$ and $q_{n + 1}$
by the case $n = 1.$
This \pj\  satisfies the conclusion of the lemma.
\end{proof}

\begin{cor}\label{CorOfZ3}
Let $A$ be a \ca\  with real rank zero.
Let $q_0, q_1, \ldots, q_n \in A$ be \pj s, and let $\ep > 0.$
Then there exists a \pj\  $e \in A$ such that:
\begin{itemize}
\item[(1)]
$e \leq q_0.$
\item[(2)]
$\| q_k e - e \| < \ep$ for $1 \leq k \leq n.$
\item[(3)]
$[e] \geq [1] - \sum_{k = 0}^n ([1] - [q_k])$ in $K_0 (A).$
\end{itemize}
\end{cor}

\begin{proof}
Apply Lemma~\ref{Z3} with the same value of $\ep,$
and with $1 - q_0, \, 1 - q_1, \, \ldots, \, 1 - q_n$
in place of $q_0, q_1, \ldots, q_n.$
Call the resulting \pj\  $f.$
Then take $e = 1 - f.$
\end{proof}

The next several lemmas will be used to control the traces of
various \pj s we construct.

\begin{lem}\label{Q2}
Let $\ep > 0.$
Then there is a \cfn\  $g \colon [0, 1] \to [0, 1]$
with $g (0) = 0$ and $g (1) = 1$ such that whenever
$A$ is a \ca\  with real rank zero,
and $a \in A$ satisfies $\| a \| \leq 1,$
then there exists a \pj\  $e \in A$ such that
$g (a a^*) e = e$ and $\| e a - a \| < \ep.$
\end{lem}

\begin{proof}
Set $\ep_0 = \ts{ \frac{1}{18}} \ep^2.$
Define \cfn s  $g, h \colon [0, 1] \to [0, 1]$ by
\[
g (t) = \left\{ \begin{array}{ll}
        \ep_0^{-1} t  & \hspace{3em} 0 \leq t \leq \ep_0  \\
        1             & \hspace{3em} \ep_0 \leq t \leq 1
                   \end{array}
                   \right.
\]
and
\[
h (t) = \left\{ \begin{array}{ll}
        0             & \hspace{3em} 0 \leq t \leq \ep_0  \\
  \ep_0^{-1} (t - \ep_0) & \hspace{3em} \ep_0 \leq t \leq 2 \ep_0 \\
        1             & \hspace{3em} 2 \ep_0 \leq t \leq 1.
                   \end{array}
                   \right.
\]
Then $\| h (a a^*) a a^* - a a^* \| < 2 \ep_0.$
Therefore Lemma~\ref{Calc} gives
\[
\| h (a a^*) a - a \| < \sqrt{2 \ep_0} = \ts{ \frac{1}{3} } \ep.
\]
Since $A$ has real rank zero,
there is a \pj\  $e$ in the \hsa\  $B$ of $A$ generated by $h (a a^*)$
such that $\| e h (a a^*) - h (a a^*) \| < \ts{ \frac{1}{3} } \ep.$
Now
\begin{align*}
\| e a - a \|
 & \leq \| e \| \cdot \| a - h (a a^*) a \|
            + \| e h (a a^*) - h (a a^*) \| \cdot \| a \|
            + \| h (a a^*) a - a \|   \\
 & < \ts{ \frac{1}{3} } \ep
            + \ts{ \frac{1}{3} } \ep + \ts{ \frac{1}{3} } \ep
   = \ep.
\end{align*}
Also $g (a a^*) e = e$ because $g (a a^*) b = b$ for all $b \in B.$
\end{proof}

\begin{lem}\label{Q1}
Let $g \colon [0, 1] \to [0, 1]$ be a \cfn\  such that
$g (0) = 0$ and $g (1) = 1,$ and let $\ep > 0.$
Then:
\begin{itemize}
\item[(1)]
There exists $\dt > 0$ such that whenever $\om$ is a state
on a \ca\  $A$ and $a \in A$ satisfies $0 \leq a \leq 1$
and $\om (a) < \dt,$ then $\om (g (a)) < \ep.$
\item[(2)]
There exists $\dt > 0$ such that whenever $\om$ is a state
on a \ca\  $A$ and $a \in A$ satisfies $0 \leq a \leq 1$
and $\om (a) > 1 - \dt,$ then $\om (g (a)) > 1 - \ep.$
\end{itemize}
\end{lem}

\begin{proof}
We prove~(1).
Choose $\dt_0 > 0$ such that $g (t) < \ts{ \frac{1}{2} } \ep$
for all $t \in [0, \dt_0].$
Set $\dt = \ts{ \frac{1}{2} } \ep \dt_0.$
Let $a$ and $\om$ be as in the hypotheses of~(1).

By considering the \hm\  $\ph \colon C ([0, 1]) \to A$
which sends the function $h (t) = t$ to $a,$
and the state $\om \circ \ph,$
we reduce to the case $A = C ([0, 1])$ and $a = h.$
Then there is a probability measure $\mu$ on $[0, 1]$
such that $\om (f) = \int_0^1 f \, d \mu$ for all $f \in C ([0, 1]).$
The assumption $\om (h) < \dt$ becomes
$\int_0^1 t \, d \mu (t) < \dt.$
{}From $t \geq \dt_0 \ch_{[\dt_0, 1]} (t)$ we get
$\dt_0 \mu ([\dt_0, 1]) < \dt,$
whence $\mu ([\dt_0, 1]) < \ts{ \frac{1}{2} } \ep.$
Since $g \leq 1,$
we now get $\int_{\dt_0}^1 g \, d \mu < \ts{ \frac{1}{2} } \ep.$
Also
\[
\int_0^{\dt_0} g \, d \mu < \ts{ \frac{1}{2} } \ep \mu ([0, \dt_0])
  \leq \ts{ \frac{1}{2} } \ep.
\]
Putting these together gives
$\om (g (a)) = \int_0^1 g \, d \mu < \ep.$

For~(2), apply~(1) to the function $h (t) = 1 - g (1 - t).$
If $\om (a) > 1 - \dt$ then $\om (1 - a) < \dt,$
so $\om (h (1 - a)) < \ep,$
whence $\om (g (a)) = 1 - \om (h (1 - a)) > 1 - \ep.$
\end{proof}

\begin{lem}\label{X1}
Let $g \colon [0, 1] \to [0, 1]$ be a \cfn\  such that
$g (0) = g (1) = 0,$ and let $\ep > 0.$
Then there exists $\dt > 0$ such that whenever $\om$ is a state
on a \ca\  $A$ and $a \in A$ satisfies $0 \leq a \leq 1$
and $\om (a - a^2) < \dt,$ then $\om (g (a)) < \ep.$
\end{lem}

\begin{proof}
Choose $\dt_0 > 0$ such that $\dt_0 < \frac{1}{2}$ and
$g (t) < \ts{ \frac{1}{3} } \ep$
for all $t \in [0, \dt_0] \cup [1 - \dt_0, \, 1].$
Set $\dt = \ts{ \frac{1}{6} } \ep \dt_0.$
Let $a$ and $\om$ be as in the hypotheses.
As in the proof of Lemma~\ref{Q1},
we may assume that $A = C ([0, 1]),$ that $a (t) = t$ for all $t,$
and that there is a probability measure $\mu$ on $[0, 1]$
such that $\om (f) = \int_0^1 f \, d \mu$ for all $f \in C ([0, 1]).$
For $t \in [\dt_0, \, 1 - \dt_0]$ we have
$t - t^2 \geq \dt_0 - \dt_0^2 > \ts{ \frac{1}{2} } \dt_0.$
{}From $\int_0^1 (t - t^2) \, d \mu (t) < \ts{ \frac{1}{6} } \ep \dt_0$
we therefore get $\mu ([\dt_0, \, 1 - \dt_0]) < \ts{ \frac{1}{3} } \ep.$
So
\begin{align*}
\om (g (a))
 & = \int_0^{\dt_0} g \, d \mu
          + \int_{\dt_0}^{1 - \dt_0} g \, d \mu
          + \int_{1 - \dt_0}^1 g \, d \mu   \\
 & \leq \ts{ \frac{1}{3} } \ep \mu ([0, \dt_0])
          + \mu ([\dt_0, \, 1 - \dt_0])
          + \ts{ \frac{1}{3} } \ep \mu ([1 - \dt_0, \, 1])
   < \ep.
\end{align*}
This completes the proof.
\end{proof}

We will make extensive use of the $L^2$-norm (or seminorm)
associated with a tracial state $\ta$ of a \ca\  $A,$ given by
$\| a \|_{2, \ta} = \ta (a^* a)^{1/2}.$
See the discussion before Lemma V.2.20 of~\cite{Tk}
for more on this seminorm in the von Neumann algebra context.
All the properties we need are immediate from its identification
with the seminorm in which one completes $A$ to obtain the
Hilbert space $H_{\ta}$ for the
Gelfand-Naimark-Segal representation associated with $\ta,$
and from the relation $\ta (b a) = \ta (a b).$
In particular, we always have
$\| a b c \|_{2, \ta} \leq \| a \| \cdot \| b \|_{2, \ta} \cdot \| c \|.$

We prove here $L^2$ analogs for \ca s with real rank zero
of various norm approximation results.
These are probably well known in von Neumann algebras.
However, as we see here, all that is really needed is real rank zero.

\begin{lem}\label{Q4}
For every $\ep > 0$ and $n \in \N,$ there is $\dt > 0$
such that the following holds.
Let $A$ be a unital \ca\  with real rank zero,
let $T \subset T (A),$
and let $p, f_1, f_2, \ldots, f_n \in A$ be \pj s
such that $\| p f_k \|_{2, \ta} < \dt$
for $1 \leq k \leq n$ and all $\ta \in T.$
Then there exists a \pj\  $e \in A$
such that:
\begin{itemize}
\item[(1)]
$e \leq p.$
\item[(2)]
$\| e f_k \| < \ep$ for $1 \leq k \leq n.$
\item[(3)]
$\ta (e) > \ta (p) - \ep$ for all $\ta \in T.$
\end{itemize}
\end{lem}

\begin{proof}
Choose $g$ as in Lemma~\ref{Q2}
for $\ts{ \frac{1}{3} } \ep$ in place of $\ep.$
Apply Lemma~\ref{Q1}(1) with this $g$
and with $n^{-1} \ep$ in place of $\ep,$
and let $\dt_0$ be the resulting positive number.
Set $\dt = \dt_0^{1/2}.$

Let the \ca\  $A,$ the subset $T \subset T (A),$
and the \pj s $p, f_1, f_2, \ldots, f_n \in A$
be as in the hypotheses, with this choice of $\dt.$
By the choice of $g,$ for each $k$ there is a \pj\  $q_k$
such that $g (p f_k p) q_k = q_k$
and $\| q_k p f_k - p f_k \| < \ts{ \frac{1}{3} } \ep.$
Apply Lemma~\ref{Z3} in the \ca\  $p A p$ to find a \pj\  $q \leq p$
such that $\| q q_k - q_k \| < \ts{ \frac{1}{3} } \ep$
for $1 \leq k \leq n$ and
$q \precsim q_1 \oplus q_2 \oplus \cdots \oplus q_n$ in $M_n (p A p).$
Then set $e = p - q.$

That $e \leq p$ is clear.
To prove~(2), we estimate,
using $q \leq p$ at the second step,
\begin{align*}
\| e f_k \|
 & = \| p f_k - q f_k \|
   = \| p f_k - q p f_k \|  \\
 & \leq \| p f_k - q_k p f_k \|
            + \| q_k - q q_k \| \cdot \| p f_k \|
            + \| q \| \cdot \| q_k p f_k - p f_k \|               \\
 & < \ts{ \frac{1}{3} } \ep
            + \ts{ \frac{1}{3} } \ep + \ts{ \frac{1}{3} } \ep
   = \ep.
\end{align*}
For~(3),
we use the choice of $\dt_0$ and the hypothesis
$\ta (p f_k p) = \| p f_k \|_{2, \ta}^2 < \dt^2 = \dt_0$
to conclude that $\ta (g (p f_k p)) < n^{-1} \ep.$
So $g (p f_k p) q_k = q_k$ gives
\[
\ta (q_k) = \ta ( g (p f_k p)^{1/2} q_k g (p f_k p)^{1/2} )
  \leq \ta (g (p f_k p)) < n^{-1} \ep.
\]
Therefore $\ta (q) \leq \sum_{k = 1}^n \ta (q_k) < \ep,$
whence $\ta (e) = \ta (p) - \ta (q) > \ta (p) - \ep.$
\end{proof}

\begin{lem}\label{X2}
For every $\ep > 0$ and $n \in \N,$ there is $\dt > 0$ such that
whenever $A$ is a \ca\  with real rank zero,
$T \subset T (A),$ $r \in A$ is a \pj,
and $p_1, p_2, \ldots, p_n \in r A r$ are \pj s
with $\| p_j p_k \|_{2, \ta} < \dt$ for $j \neq k$ and $\ta \in T,$
then there exist \mops\  $q_1, q_2, \ldots, q_n \in r A r$
such that
$\| q_k - p_k \|_{2, \ta} < \ep$ for all $k$ and all $\ta \in T.$
\end{lem}

\begin{proof}
We prove this by induction on $n.$
The statement is trivial for $n = 1.$
Suppose therefore that the statement is known for a
particular value of $n,$ and let $\ep > 0.$
Choose $\dt_0$ for the statement with $n$ \pj s
and with $\ts{ \frac{1}{4}} \ep$ in place of $\ep.$
Set
$\dt_1 = \min \left( \ts{ \frac{1}{3}} \dt_0,
                 \, \ts{ \frac{1}{2}} \ep \right).$
Choose $\dt_2 > 0$ so small that whenever
$f, g$ are \pj s in a \ca\  $B$ which satisfy
$\| f g \| < \dt_2,$
then there is a \pj\  $h \in B$
such that $f h = 0$ and $\| g - h \| < \dt_1.$
We also require $\dt_2 \leq \ep.$
Choose $\dt > 0$ as in Lemma~\ref{Q4} for the given value of $n$
and with $\dt_2$ in place of $\ep.$
We also require $\dt \leq \ts{ \frac{1}{3}} \dt_0.$

Now let $A,$ $T,$ and $r, \, p_1, \, p_2, \, \ldots, \, p_{n + 1} \in A$
be as in the hypotheses, with this value of $\dt.$
Use the choice of $\dt$ following Lemma~\ref{Q4}
to find a \pj\  $q_{n + 1} \in A$ with
$q_{n + 1} \leq p_{n + 1}$ such that
$\| q_{n + 1} p_k \| < \dt_2$ for $1 \leq k \leq n$ and
$\ta (p_{n + 1} - q_{n + 1}) < \dt_2.$
Then also $\| [q_{n + 1} + (1 - r)] p_k \| < \dt_2$
for $1 \leq k \leq n.$
Use the choice of $\dt_2$
to find \pj s $e_1, e_2, \ldots, e_n \in A$
which are orthogonal to $q_{n + 1} + (1 - r)$
and satisfy $\| e_k - p_k \| < \dt_1$ for $1 \leq k \leq n.$

For $j \neq k$ we now have
\[
\| e_j e_k \|_{2, \ta}
 \leq \| e_j - p_j \| \cdot \| e_k \|_{2, \ta}
        + \| p_j \|_{2, \ta} \cdot \| e_k - p_k \|
        + \| p_j p_k \|_{2, \ta}
 < 2 \dt_1 + \dt
 \leq \dt_0.
\]
Apply the induction hypothesis with $e_1, e_2, \ldots, e_n$
in place of $p_1, p_2, \ldots, p_n$ and with $r - q_{n + 1}$
in place of $r,$
obtaining \mops\  $q_1, q_2, \ldots, q_n \leq r - q_{n + 1}$
with $\| q_k - e_k \|_{2, \ta} < \ts{ \frac{1}{4}} \ep$
for all $k$ and all $\ta \in T.$
Then for $1 \leq k \leq n$ and $\ta \in T$ we have
\[
\| q_k - p_k \|_{2, \ta}
   \leq \| q_k - e_k \|_{2, \ta} + \| e_k - p_k \|
   < \ts{ \frac{1}{4}} \ep + \dt_1
   \leq \ep.
\]
Since also $\| q_{n + 1} - p_{n + 1} \|_{2, \ta} < \dt_2 \leq \ep,$
this completes the proof.
\end{proof}

The next two lemmas are special cases of Lemma~\ref{X5},
but are used in its proof.
In the first, we could presumably, with more care,
also require $\| p q - q \| < \ep.$

\begin{lem}\label{X3}
For every $\ep > 0$ and $n \in \N,$ there is $\dt > 0$ such that
whenever $A$ is a unital \ca\  with real rank zero,
$T \subset T (A),$ and $(e_{j, k})_{1 \leq j, k \leq n}$
is a system of matrix units for a copy of $M_n$ in $A$
with identity $e,$
and whenever $p \in  A$ is a \pj\  with $p \leq e$
such that $\| [p, e_{j, k}] \|_{2, \ta} < \dt$
for $1 \leq j, k \leq n$ and all $\ta \in T,$
then there is a \pj\  $q \in A$ with $q \leq e$
which commutes with all $e_{j, k}$
and such that $\| q - p \|_{2, \ta} < \ep$ for all $\ta \in T.$
\end{lem}

\begin{proof}
Use Lemma~5.1 of~\cite{OP1}
to choose a \cfn\  $g \colon [0, 1] \to [0, 1]$
such that $g (0) = 0,$ such that $g (1) = 1,$ and such that
whenever $A$ is a \ca\  with real rank zero and $a \in A$
satisfies $0 \leq a \leq 1,$
then there is a \pj\  $f \in {\overline{a A a}}$ such that
\[
\| a f - f \| < \frac{\ep}{6 n} \andeqn
\| f g (a) - g (a) \| < \frac{\ep}{6 n}.
\]
Apply Lemma~\ref{X1} to the function
$h (t) = (t - g (t))^2$ with $\frac{1}{144} n^{- 2} \ep^2$
in place of $\ep.$
Let $\dt_0 > 0$ be the resulting number.
Set
\[
\dt = \min \left( \sqrt{2 \dt_0}, \, \frac{\ep}{2 n^2} \right).
\]

Now let $A,$ $T,$ $(e_{j, k})_{1 \leq j, k \leq n},$ $e,$ and $p$
be as in the statement of the lemma.
Then, using the trace property several times at the last step,
\begin{align*}
2 \dt_0 & \geq \dt^2
   > \ta ([p, e_{1, 1}]^* [p, e_{1, 1}])  \\
 & = \ta ( e_{1, 1} p e_{1, 1} + p e_{1, 1} p
             - e_{1, 1} p e_{1, 1} p - p e_{1, 1} p e_{1, 1} )   \\
 & = 2 \ta (e_{1, 1} p e_{1, 1} - ( e_{1, 1} p e_{1, 1})^2 ).
\end{align*}
With $a = e_{1, 1} p e_{1, 1},$ the choice of $\dt_0$ gives
$\ta ((g (a) - a)^2) < \frac{1}{144} n^{- 2} \ep^2$
for all $\ta \in T.$
The choice of $g$ then gives a \pj\  %
$f \in {\overline{a A a}} \subset e_{1, 1} A e_{1, 1}$ such that
\[
\| a f - f \| < \frac{\ep}{6 n} \andeqn
\| f g (a) - g (a) \| < \frac{\ep}{6 n}.
\]
We estimate:
\begin{align*}
\| f - a \|_{2, \ta}
 & \leq \| f - f a \| + \| f \| \cdot \| a - g (a) \|_{2, \ta}
        + \| f g (a) - g (a) \|  + \| g (a) - a \|_{2, \ta}      \\
 & < \frac{\ep}{6 n}
        + \left( \frac{\ep^2}{144 n^2} \right)^{1/2}
        + \frac{\ep}{6 n}
        + \left( \frac{\ep^2}{144 n^2} \right)^{1/2}
   = \frac{\ep}{2 n}.
\end{align*}

Now set $q = \sum_{k = 1}^n e_{k, 1} f e_{1, k}.$
Clearly $q$ commutes with all $e_{j, k}$
and $q \leq e = \sum_{k = 1}^n e_{k, k}.$
We estimate $\| p - q \|_{2, \ta},$ in several steps.
First, for $j \neq k$ we have
\[
\| e_{j, j} p e_{k, k} \|_{2, \ta}
  \leq \| [e_{j, j}, p] \|_{2, \ta} \| e_{k, k} \|
  < \dt \leq \frac{\ep}{2 n^2}.
\]
Second,
we get
\begin{align*}
\| e_{k, k} p e_{k, k} - e_{k, 1} f e_{1, k} \|_{2, \ta}
 & \leq \| e_{k, k} p e_{k, k} - e_{k, 1} p e_{1, k} \|_{2, \ta}
          + \| e_{1, 1} p e_{1, 1} - f \|_{2, \ta}              \\
 & \leq \| e_{k, k} \| \cdot \| [p, e_{k, 1}] \|_{2, \ta}
                \cdot \| e_{1, k} \|
          + \| f - a \|_{2, \ta}    \\
 & < \dt + \frac{\ep}{2 n}
   \leq \frac{\ep}{2 n^2} + \frac{\ep}{2 n}.
\end{align*}
Adding things up, we get
\begin{align*}
\| p - q \|_{2, \ta}
 & \leq \sum_{k = 1}^n
         \| e_{k, k} p e_{k, k} - e_{k, 1} f e_{1, k} \|_{2, \ta}
     + \sum_{j \neq k} \| e_{j, j} p e_{k, k} \|_{2, \ta}         \\
 & < n \left( \frac{\ep}{2 n^2} + \frac{\ep}{2 n} \right)
        + (n^2 - n) \left( \frac{\ep}{2 n^2} \right)
   = \ep.
\end{align*}
This completes the proof.
\end{proof}

\begin{lem}\label{X4}
For every $\ep > 0$ and $n \in \N,$ there is $\dt > 0$ such that
whenever $A$ is a unital \ca\  with real rank zero,
$T \subset T (A),$ and $e_1, e_2, \ldots, e_n \in A$
are \pj s such that $\sum_{k = 1}^n e_k = 1,$
and whenever $p \in  A$ is a \pj\  %
such that $\| [p, e_k] \|_{2, \ta} < \dt$
for $1 \leq j, k \leq n$ and all $\ta \in T,$
then there is a \pj\  $q \in A$ which commutes with all $e_k$
and such that $\| q - p \|_{2, \ta} < \ep$ for all $\ta \in T.$
\end{lem}

\begin{proof}
The proof is essentially the same as that of Lemma~\ref{X3}.
Let $g,$ $\dt_0,$ and $\dt$ be as there.
With $A,$ $T,$ $e_1, e_2, \ldots, e_n,$ and $p$ as in
the hypotheses,
for $1 \leq k \leq n$ let $a_k = e_k p e_k,$
and use the method of the first part of the proof of Lemma~\ref{X3}
to find \pj s $f_k \in e_k A e_k$ such that
$\| f_k - a_k \|_{2, \ta} < \frac{1}{2} n^{-1} \ep.$
Set $q = \sum_{k = 1}^n f_k.$
Using the method of the second part of the proof of Lemma~\ref{X3},
slightly simplified,
one gets
\[
\| q - p \|_{2, \ta}
 < n \left(\frac{\ep}{2 n} \right)
        + (n^2 - n) \left( \frac{\ep}{2 n^2} \right)
   < \ep,
\]
as desired.
\end{proof}

\begin{lem}\label{X5}
For every $\ep > 0$ and \fd\  \ca\  $E,$ there is $\dt > 0$ such that
whenever $A$ is a unital \ca\  with real rank zero,
$T \subset T (A),$ $S$ is a system of matrix units for $E,$
and $\ph \colon E \to A$ is a unital \hm,
and whenever $p \in  A$ is a \pj\  %
such that $\| [p, \, \ph (x)] \|_{2, \ta} < \dt$
for all $x \in S$ and $\ta \in T,$
then there is a \pj\  $q \in A$ which commutes with $\ph (E)$
and such that $\| q - p \|_{2, \ta} < \ep$ for all $\ta \in T.$
\end{lem}

\begin{proof}
Write $E = \bigoplus_{k = 1}^n M_{r (k)}.$
Apply Lemma~\ref{X3}
with $\ts{ \frac{1}{2}} n^{-1} \ep$ in place of $\ep$
and with the values
$r (1), \, r (2), \, \ldots, \, r (n)$ in place of $n.$
Let $\dt_0 > 0$ be the minimum of $\ep$
and the $n$ resulting values of $\dt.$
Apply Lemma~\ref{X4} with $\ts{ \frac{1}{3}} \dt_0$ in place of $\ep$
and with the given value of $n,$
and let $\dt_1 > 0$ be the resulting number.
Set $R = \max (r (1), \, r (2), \, \ldots, \, r (n) ),$
and set
$\dt = \min \left( R^{-1} \dt_1, \, \ts{ \frac{1}{3}} \dt_0 \right).$

Now let $A,$ $T,$ $E,$ $S,$ and $p$
be as in the statement of the lemma.
Let $e_1, e_2, \ldots, e_n$ be the minimal central \pj s of $E,$
with $e_k$ being the identity of $M_{r (k)}.$
Since $e_k$ is a sum of $r (k) \leq R$ elements of $S,$
we get $\| [p, \, \ph (e_k)] \|_{2, \ta} < R \dt \leq \dt_1$
for $1 \leq k \leq n$ and for all $\ta \in T.$
Therefore there is a \pj\  $f \in A$ which commutes with all
$\ph (e_k)$ and such that
$\| f - p \|_{2, \ta} < \ts{ \frac{1}{3}} \dt_0$ for all $\ta \in T.$
Set $f_k = f \ph (e_k).$
For $s \in S$ a matrix unit for the $k$-th summand $M_{r (k)},$
we have
$[\ph (s), f_k] = [\ph (s), f],$
so
\[
\| [\ph (s), f_k] \|_{2, \ta}
 \leq 2 \| f - p \|_{2, \ta} \| \ph (s) \| + \| [\ph (s), p] \|_{2, \ta}
 < \ts{ \frac{2}{3}} \dt_0 + \ts{ \frac{1}{3}} \dt_0
 = \dt_0.
\]
Applying Lemma~\ref{X3} to the set of all such $\ph (s)$
with $k$ fixed, and with $e_k$ in place of $e,$
we obtain a \pj\  $q_k \leq e_k$ which commutes with all
such $\ph (s)$ and satisfies
$\| q_k - f_k \|_{2, \ta} < \ts{ \frac{1}{2}} n^{-1} \ep.$
Set $q = \sum_{k = 1}^n q_k,$
which is a \pj\  which commutes with all elements of $E$
and satisfies
\[
\| q - f \|_{2, \ta}
 \leq \sum_{k = 1}^n \| q_k - f_k \|_{2, \ta}
 < n \left( \ts{ \frac{1}{2}} n^{-1} \ep \right)
 = \ts{ \frac{1}{2}} \ep.
\]
So
$\| q - p \|_{2, \ta}
 < \ts{ \frac{1}{2}} \ep + \ts{ \frac{1}{3}} \dt_0
 \leq \ep.$
\end{proof}

\begin{lem}\label{Z2pp}
Let $A$ be a simple separable unital \ca\  with tracial rank zero.
Let $\af \colon \Gm \to \Aut (A)$  be an action of a discrete group.
Let $F \subset \Gm$ and $S \subset A$ be finite sets,
and let $\ep > 0.$
Then there exist \pj s $q, q_0 \in A,$
unital \fd\  subalgebras $E \subset q A q$ and $E_0 \subset q_0 A q_0,$
and automorphisms $\ph_{\gm} \in \Aut (A)$
for $\gm \in F \cup \{ 1 \},$ such that:
\begin{itemize}
\item[(1)]
$\ph_1 = \id_A$ and
$\| \ph_{\gm} - \af_{\gm} \| < \ep$ for all $\gm \in F.$
\item[(2)]
For every  $\gm \in F \cup \{ 1 \}$ and $x \in E,$
we have $q_0 \ph_{\gm} (x) = \ph_{\gm} (x) q_0$
and $q_0 \ph_{\gm} (x) q_0 \in E_0.$
\item[(3)]
For every $a \in S,$ we have $\| q a - a q \| < \ep$ and
$\dist (q a q, \, E) < \ep.$
\item[(4)]
$\ta (1 - q), \, \ta (1 - q_0) < \ep$ for all $\ta \in T (A).$
\end{itemize}
\end{lem}

\begin{proof}
Since $A$ has tracial rank zero, there is a \pj\  $q \in A$
and a unital \fd\  subalgebra $E \subset q A q,$
such that $\ta (1 - q) < \ep$ for all $\ta \in T (A),$
and such that for every $a \in S,$
we have $\| q a - a q \| < \ep$ and $\dist (q a q, \, E) < \ep.$
With this choice, Condition~(3) clearly holds,
as does the part of Condition~(4) which applies to $q.$

Use Lemma~2.5.10 of~\cite{LnBook} to find $\dt > 0$
such that whenever $B$ and $C$ are unital subalgebras
of a unital \ca\  $A$ such that $\dim (B) \leq \dim (E) + 1$
and such that $B$ has a system of matrix units each of
which has distance less than $\dt$ from $C,$ then there is
a unitary $u \in A$ with $\| u - 1 \| < \ts{ \frac{1}{4}} \ep$
such that $u^* B u \subset C.$

Let $T_0$ be a system of matrix units for $E + \C (1 - q),$
and set $T = \bigcup_{\gm \in F \cup \{ 1 \}} \af_{\gm} (T_0).$
Using again the fact that $A$ has tracial rank zero,
find a \pj\  $p \in A$
and a unital \fd\  subalgebra $D \subset p A p$
such that $\ta (1 - p) < \ep$ for all $\ta \in T (A),$
and such that for every $a \in T,$
we have $\| p a - a p \| < \ts{ \frac{1}{3}} \dt$
and $\dist (p a p, \, D) < \ts{ \frac{1}{3}} \dt.$
Set $B = D + (1 - p) A (1 - p),$ and observe that if $a \in T$
then
\begin{align*}
\dist (a, B)
 & \leq \| a - [p a p + (1 - p) a (1 - p)] \| + \dist (p a p, \, D)
              \\
 & \leq 2 \| [p, a] \| + \dist (p a p, \, D) < \dt.
\end{align*}
Using the choice of $\dt$ on the system of matrix units
$\af_{\gm} (T_0)$ for the subalgebra $\af_{\gm} (E + \C (1 - q)),$
choose for each $\gm \in F \cup \{ 1 \}$ a unitary $u_{\gm} \in A$
with $\| u_{\gm} - 1 \| < \ts{ \frac{1}{4}} \ep$
such that
\[
u_{\gm}^* \af_{\gm} (E + \C (1 - q)) u_{\gm}
   \subset D + (1 - p) A (1 - p).
\]
Now set $q_0 = u_1 p u_1^*$ and $E_0 = u_1 D u_1^*.$
This gives
\[
E + \C (1 - q) \subset E_0 + (1 - q_0) A (1 - q_0).
\]
Further define
$\ph_{\gm} (x) = [u_1 u_{\gm}^*] \af_{\gm} (x) [u_1 u_{\gm}^*]^*$
for $x \in E.$
Then $\ph_1 = \id_A$ and
\[
\| \ph_{\gm} - \af_{\gm} \|
 \leq 2 \| u_1 - 1 \| + 2 \| u_{\gm} - 1 \| < \ep.
\]
Thus Condition~(1) holds.
Condition~(2) follows from
\[
u_1 [D + (1 - p) A (1 - p)] u_1^* = E_0 + (1 - q_0) A (1 - q_0),
\]
and the remaining part of Condition~(4) holds by construction.
\end{proof}

In following theorem, Condition~(3) can also be weakened
to require merely $\sum_{k = 0}^n \ta (e_k) > 1 - \ep$
for every $\ta \in T (A).$

\begin{thm}\label{TraceVersion}
Let $A$ be a simple separable unital \ca\  with tracial rank zero.
Let $\af \in \Aut (A).$
Then $\af$ generates an action of $\Z$ with the \tRp\  \ifo\  %
for every $\ep > 0,$ every $n \in \N,$
and every finite subset $S \subset A$ there exist
orthogonal \pj s $e_0, e_1, \ldots, e_n \in A$ such that:
\begin{itemize}
\item[(1)]
$\| \af (e_k) - e_{k + 1} \|_{2, \ta} < \ep$
for $0 \leq k \leq n - 1$ and all $\ta \in T (A).$
\item[(2)]
$\| [e_k, a] \|_{2, \ta} < \ep$
for $0 \leq k \leq n,$ all $a \in S,$ and all $\ta \in T (A).$
\item[(3)]
$\sum_{k = 0}^n e_k = 1.$
\end{itemize}
\end{thm}

\begin{proof}
Suppose $\af$ generates an action of $\Z$ with the \tRp.
Let $\ep,$ $S,$ and $n$ be given.
We may assume that $\| a \| \leq 1$ for all $a \in S.$
Since $A$ has tracial rank zero,
Corollary~5.7 and Theorems~5.8 and~6.8 of~\cite{LnTTR}
imply that the order on \pj s over $A$ is determined by traces.
We apply Lemma~1.4 of~\cite{OP1} with
$\min \left( \ts{ \frac{1}{3}} \ep, \ts{ \frac{1}{9}} \ep^2 \right)$
in place of $\ep.$
Call the resulting \pj s $f_0, f_1, \ldots, f_n.$
Set $f = \sum_{j = 0}^n f_j,$ set $e_0 = f_0 + 1 - f,$
and set $e_j = f_j$ for $1 \leq j \leq n.$
Since
$\| e_0 - f_0 \|_{2, \ta} = \ta (1 - f)^{1/2} < \ts{ \frac{1}{3}} \ep$
and $\| a \|_{2, \ta} \leq \| a \|$
for all $\ta \in T (A),$ it is easy to show that
$e_0, e_1, \ldots, e_n$ satisfy the conditions of the theorem.

Now we prove the other direction.
Let $S \subset A$ be a finite set, let $\ep > 0,$ and let $n \in \N.$
\Wolog\  $\ep < 1.$
We construct \pj s $e_0, e_1, \ldots, e_n$ such that
$\| \af (e_k) - e_{k + 1} \| < \ep$ for $0 \leq k \leq n - 1,$
such that $\| [e_k, a] \| < \ep$ for $a \in S$ and $0 \leq k \leq n,$
and such that
$1 - \sum_{k = 0}^n \ta (e_k) < \ep$ for all $\ta \in T (A).$
By Lemma~1.4 of~\cite{OP1}, this is sufficient.

Choose $\ep_0 > 0$ with $\ep_0 < \ts{ \frac{1}{3}} \ep (n + 1)^{-1}$
and so small that whenever
$p_0, p_1, \ldots, p_n$ are \pj s in a \ca\  $B$ which satisfy
$\| p_j p_k \| < 4 \ep_0$ for $j \neq k,$
then there are orthogonal \pj s $e_0, e_1, \ldots, e_n \in B$
such that $e_0 = p_0$ and $\| e_k - p_k \| < \ts{ \frac{1}{3}} \ep$
for $1 \leq k \leq n.$

Apply Lemma~\ref{Q4} with $\ep_0$ in place of $\ep$ and with $n$ as
given, and let $\ep_1 > 0$ be the resulting value of $\dt.$
We also require $\ep_1 \leq \ep_0.$
Then set
\[
\ep_2 = \min \left(1, \, \frac{\ep_1^2}{16}, \, \frac{\ep_1}{2 (n + 3)},
      \, \frac{\ep}{18} \right).
\]
Set
\[
T = \bigcup_{k = 0}^n \af^{-k} (S).
\]
Apply Lemma~\ref{Z2pp} with $\ep_2$ in place of $\ep,$
with $T$ in place of $S,$
and with $\{ -n, \, - n + 1, \, \ldots, \, 1, \, 0 \}$
in place of $F.$
We obtain \pj s $q, q_0 \in A,$
unital \fd\  subalgebras $E \subset q A q$ and $E_0 \subset q_0 a q_0,$
and automorphisms $\ps_k = \ph_{- k} \in \Aut (A)$
for $- n \leq k \leq 0$
as there.
Apply Lemma~\ref{X5} with $\ep_2$
in place of $\ep$ and $E_0 + \C (1 - q_0)$ in place of $E,$
obtaining $\dt > 0.$

Apply the hypothesis with $n$ as given,
with $\dt$ in place of $\ep,$
and with a system of matrix units for $E_0 + \C (1 - q_0)$
in place of $S,$
getting \pj s $p_0, p_1, \ldots, p_n.$
Let $B = A \cap [E_0 + \C (1 - q_0)]',$
the subalgebra of $A$ consisting of all elements which commute
with everything in $E_0 + \C (1 - q_0).$
Apply the choice of $\dt$ using Lemma~\ref{X5} to $p_0,$
obtaining a \pj\  $f \in B$
which satisfies $\| f - p_0 \|_{2, \ta} < \ep_2$ for all $\ta \in T (A).$
Since $q_0$ is in the center of $E_0 + \C (1 - q_0),$
the element $f_0 = q_0 f$ is also a \pj\  in $B.$
Since $q_0 \ps_k (E) q_0 \subset E_0,$
it follows that $f_0$ commutes with all elements of $q_0 \ps_k (E) q_0$
and hence with all elements of $\ps_k (E).$
Therefore $f_k = \ps_k^{-1} (f_0)$ commutes with all elements of $E,$
including $q.$
So $f_k$ also commutes with $1 - q.$
We estimate $\| f_0 f_k \|_{2, \ta}$ for $\ta \in T (A).$
To start, $\| f - p_0 \|_{2, \ta} < \ep_2$ and
\[
\| f_0 - f \|_{2, \ta} = \ta ( (1 - q_0) f)^{1/2}
  \leq \ta (1 - q_0)^{1/2}
  < \ep_2^{1/2}
  \leq \ts{ \frac{1}{4}} \ep_1,
\]
so
\[
\| f_0 - p_0 \|_{2, \ta}
  \leq \| f_0 - f \|_{2, \ta} + \| f - p_0 \|_{2, \ta}
  < \ts{ \frac{1}{4}} \ep_1 + \ep_2.
\]
Now for $1 \leq k \leq n$ we get
\begin{align*}
\| f_k - \af^k (f_0) \|
& = \| \ps_k^{-1} (f_0) - \af^k (f_0) \|
  \leq \| \ps_k^{-1} - \af^k \|     \\
& \leq \| \ps_k^{-1} \| \cdot \| \af^{- k} - \ps_k \| \cdot \| \af^k \|
  < \ep_2.
\end{align*}
So
\begin{align*}
\| f_k - p_k \|_{2, \ta}
 & \leq \| f_k - \af^k (f_0) \|_{2, \ta}
           + \| \af^k (f_0) - \af^k (p_0) \|_{2, \ta}
           + \| \af^k (p_0) - p_k \|_{2, \ta}   \\
 & \leq \| f_k - \af^k (f_0) \| + \| f_0 - p_0 \|_{2, \, \ta \circ \af^k}
           + \sum_{j = 0}^{k - 1}
                  \| \af (p_j) - p_{j + 1} \|_{2, \, \ta \circ \af^j}  \\
 & < \ep_2 + \left( \ts{ \frac{1}{4}} \ep_1 + \ep_2 \right) + n \ep_2
   = \ts{ \frac{1}{4}} \ep_1 + (n + 2) \ep_2.
\end{align*}
Therefore,
\begin{align*}
\| f_0 f_k \|_{2, \ta}
 & = \| f_0 f_k - p_0 p_k \|_{2, \ta}
   \leq \| f_0 - p_0 \|_{2, \ta} \cdot \| f_k \|
           + \| p_0 \| \cdot \| f_k -p_k \|_{2, \ta}       \\
 & < \left( \ts{ \frac{1}{4}} \ep_1 + \ep_2 \right)
           + \left( \ts{ \frac{1}{4}} \ep_1 + (n + 2) \ep_2 \right)
   = \ts{ \frac{1}{2}} \ep_1 + (n + 3) \ep_2
   \leq \ep_1.
\end{align*}

We saw that $f_0, f_1, \ldots, f_n \in B = A \cap [E + \C (1 - q)]'.$
This algebra has real rank zero because $E + \C (1 - q)$ is \fd.
Therefore Lemma~\ref{Q4},
applied to $B$ with $T = \{ \ta |_B \colon \ta \in T (A) \},$
and the choice of $\ep_1,$ provide
a \pj\  $g \in A \cap [E + \C (1 - q)]'$
such that $g \leq f_0,$
such that $\| g f_k \| < \ep_0$ for $1 \leq k \leq n,$
and such that $\ta (g) > \ta (f_0) - \ep_0$
for all $\ta \in T (A).$

Now use Corollary~\ref{CorOfZ3} to find a \pj\  %
$e_0 \in A \cap [E + \C (1 - q)]'$ such that
$e_0 \leq q,$
such that $\| g e_0 - e_0 \| < \ep_0,$
and such that $[e_0] \geq 1 - (1 - [q]) - (1 - [g])$
in $K_0 (A \cap [E + \C (1 - q)]').$
The last inequality implies that
\[
\ta (e_0) \geq \ta (g) - \ta (1 - q)
  > \ta (f_0) - \ep_0 - \ep_2
  \geq \ta (f_0) - 2 \ep_0
\]
for all $\ta \in T (A).$

We now show that $e_0, \, \af (e_0), \, \ldots, \, \af^n (e_0)$
are approximately orthogonal.
It suffices to estimate $\| e_0 \af^k (e_0) \|$
for $1 \leq k \leq n.$
We first use $\ps_k^{-1} (g) \leq \ps_k^{-1} (f_0) = f_k$ to estimate
$\| \ps_k^{-1} (g) g \| \leq \| f_k g \| < \ep_0.$
Then
\[
\| \ps_k^{-1} (e_0) e_0 \|
  \leq 2 \| e_0 - g e_0 \|
        + \| \ps_k^{-1} (e_0) \| \cdot
               \| \ps_k^{-1} (g) g \| \cdot \| e_0 \|
  < 2 \ep_0 + \ep_0 = 3 \ep_0.
\]
So
\[
\| e_0 \af^k (e_0) \|
 \leq \| \af^{-k} - \ps_k \| + \| \ps_k^{-1} (e_0) e_0 \|
 < \ep_2 + 3 \ep_0 \leq 4 \ep_0.
\]
Now use the choice of $\ep_0$ to find orthogonal \pj s
\[
e_1, \ldots, e_n \in A \cap [E + \C (1 - q)]',
\]
all orthogonal to $e_0,$ such that
$\| e_k - \af^k (e_0) \| < \ts{ \frac{1}{3}} \ep$ for $1 \leq k \leq n.$
We prove that $e_0, e_1, \ldots, e_n$ satisfy
the estimates asked for at the beginning of the proof
of this direction.

The first estimate is easy:
\[
\| \af (e_k) - e_{k + 1} \|
 \leq \| e_k - \af^k (e_0) \| + \| \af^{k + 1} (e_0) - e_{k + 1} \|
 < \ts{ \frac{1}{3}} \ep + \ts{ \frac{1}{3}} \ep
 < \ep.
\]

For the second, we start by estimating $\| [e_0, a] \|$ for $a \in T.$
By construction, there is $b_0 \in E$ such that
$\| q a q - b_0 \| < \ep_2.$
Then $b = b_0 + (1 - q) a (1 - q)$ satisfies
\[
\| a - b \|
 \leq \| q a q - b_0 \| + 2 \| [q, a ] \|
 < 3 \ep_2.
\]
Now $e_0 \leq q$ and $e_0$ commutes with all elements
of $E,$ so $e_0$ commutes with $b.$
Therefore $\| [ e_0, a ] \| < 2 \cdot 3 \ep_2.$
If now $a \in S$ and $0 \leq k \leq n,$ then
$\af^{-k} (a) \in T,$ so
\[
\| [ \af^k (e_0), \, a ] \| = \| [ e_0, \, \af^{-k} (a) ] \|
 < 6 \ep_2,
\]
and
\[
\| [e_k, a] \| \leq 2 \| e_k - \af^k (e_0) \| + \| [ \af^k (e_0), \, a ] \|
 < \ts{ \frac{2}{3}} \ep + 6 \ep_2 \leq \ep,
\]
as desired.

Finally, we estimate $1 - \sum_{k = 0}^n \ta (e_k)$
for $\ta \in T (A).$
We saw above that $\ta (e_0) > \ta (f_0) - 2 \ep_0$
for all $\ta \in T (A).$
Thus also $\ta (\af^k (e_0)) > \ta (\af^k (f_0)) - 2 \ep_0$
for $\ta \in T (A)$ and $0 \leq k \leq n.$
{}From $\| e_k - \af^k (e_0) \| < \frac{1}{3} \ep < \frac{1}{3}$
we get $\ta (e_k) = \ta (\af^k (e_0)),$
and from $\| f_k - \af^k (f_0) \| < \ep_2 \leq 1$
we get $\ta (f_k) = \ta (\af^k (f_0)).$
So $\ta (e_k) > \ta (f_k) - 2 \ep_0$
for $\ta \in T (A)$ and $0 \leq k \leq n.$

Next, using a previous estimate at the second step,
\[
| \ta (f_k) - \ta (p_k) | \leq \| f_k - p_k \|_2
 < \ts{ \frac{1}{4}} \ep_1 + (n + 2) \ep_2 \leq \ep_1,
\]
so $\ta (f_k) > \ta (p_k) - \ep_1.$
Thus $\ta (e_k) > \ta (p_k) - \ep_1 - 2 \ep_0 \geq \ta (p_k) - 3 \ep_0.$
Summing over $k$ and using $\sum_{k = 0}^n \ta (p_k) = 1,$
we get $\sum_{k = 0}^n \ta (e_k) > 1 - 3 (n + 1) \ep_0 \geq 1 - \ep,$
as desired.
\end{proof}

\begin{lem}\label{Z1}
Let $H$ be a Hilbert space, and let $A \subset L (H)$ be a unital
C*-subalgebra which has real rank zero.
Let $p \in A''$ be a \pj.
Then $p$ is a strong operator limit of \pj s in $A.$
\end{lem}

\begin{proof}
The proof of Lemma~4.6 of~\cite{Ks2}
(which gives this conclusion for the weak operator closure
of an AF~algebras in the Gelfand-Naimark-Segal representation
from a tracial state) applies here with no essential change.
\end{proof}

The following lemma is known,
and is a special case of results in~\cite{Bn}.
However, a direct proof is much easier,
and we have not found one in the literature.

\begin{lem}\label{Hyp}
Let $A$ be a simple separable unital \ca\  with tracial rank zero,
and suppose that $A$ has a unique tracial state $\ta.$
Let $\pi_{\ta} \colon A \to L (H_{\ta})$
be the Gelfand-Naimark-Segal representation associated with $\ta.$
Then $\pi_{\ta} (A)''$ is a hyperfinite factor.
\end{lem}

\begin{proof}
Set $N = \pi_{\ta} (A)'',$ and regard $A$ as a subalgebra of $N.$
That $N$ is a factor is well known.
For hyperfiniteness, we verify~(iii) of Theorem III.7.3 of~\cite{Dx}.
The first part is trivial.
For the second part
(approximation of finite sets in trace norm by
finite dimensional subalgebras), let $t_1, t_2, \ldots, t_n \in N.$
\Wolog\  $\| t_j \| \leq 1$ for all $j.$
Let $\ep > 0.$
The sets
\[
\left\{ x \in N \colon {\mbox{$ \| x \| \leq 1$
    and $\| x - t_j \|_{2, \ta} < \tfrac{1}{2} \ep$}} \right\}
\]
are open in the *-strong operator topology on the
closed unit ball of $N$ because $\ta$ is normal
(Proposition V.2.5 of~\cite{Tk}).
So the Kaplansky Density Theorem provides
$b_1, b_2, \ldots, b_n \in A$ with
$\| b_j - t_j \|_{2, \ta} < \tfrac{1}{2} \ep$ for all $j.$
Use tracial rank zero to find a \pj\  $p \in A,$
a finite dimensional unital subalgebra $E \subset p A p,$
and $c_1, c_2, \ldots, c_n \in E$ such that
$\ta (1 - p) < \tfrac{1}{36} \ep^2,$
such that $\| p b_j - b_j p \| < \tfrac{1}{6} \ep$ for all $j,$
and such that $\| p b_j p - c_j \| < \tfrac{1}{6} \ep$ for all $j.$
Now
\[
\| b_j - p b_j p \|_{2, \ta}
  \leq \| b_j - p b_j \|_{2, \ta}
         + \| p \| \cdot \| b_j - b_j p \|_{2, \ta}
  \leq 2 \| 1 - p \|_{2, \ta}
  < \tfrac{1}{3} \ep,
\]
so
\[
\| c_j - t_j \|_{2, \ta}
  \leq \| c_j - p b_j p \| + \| p b_j p - b_j \|_{2, \ta}
              + \| b_j - t_j \|_{2, \ta}
  < \tfrac{1}{6} \ep + \tfrac{1}{3} \ep + \tfrac{1}{2} \ep
  = \ep.
\]
This completes the proof.
\end{proof}

\begin{thm}\label{OuterImpTRP}
Let $A$ be a simple separable unital \ca\  with tracial rank zero,
and suppose that $A$ has a unique tracial state $\ta.$
Let $\pi_{\ta} \colon A \to L (H_{\ta})$
be the Gelfand-Naimark-Segal representation associated with $\ta.$
Let $\af \in \Aut (A).$
Then $\af$ generates an action of $A$ with the \tRp\  \ifo\  %
for every $n > 0$ the
automorphism of $\pi_{\ta} (A)''$ induced by $\af^n$ is outer.
\end{thm}

\begin{proof}
Assume the
automorphism of $\pi_{\ta} (A)''$ induced by $\af^n$ is outer
for every $n > 0.$
We verify the hypotheses of Theorem~\ref{TraceVersion}.
Thus let $\ep > 0,$ let $n \in \N,$
and let $S \subset A$ be a finite subset.
\Wolog\  $\| a \| \leq 1$ for all $a \in S.$
Set $\ep_0 = (4 n + 5)^{-1} \ep.$
Choose $\dt > 0$ as in Lemma~\ref{X2} with $n$ as given
and for $\ep_0$ in place of $\ep.$

We regard $A$ as a subalgebra of $N = \pi_{\ta} (A)'',$
we let $\ta$ also denote the extension
of the tracial state to $\pi_{\ta} (A)'',$
and we let ${\overline{\af}}$ denote the extension
of the automorphism to $\pi_{\ta} (A)''.$
The algebra $N$ is hyperfinite by Lemma~\ref{Hyp}.

Fix $\om \in \bt \N \setminus \N,$
and let $N_{\om}$ be the central sequence algebra,
as defined before Theorem XIV.4.6 of~\cite{Tk3}.
(We use~\cite{Tk3} as our standard reference,
but almost everything is also in~\cite{Cn2}.)
By Lemma XIV.4.5 and Theorems XIV.4.6 and XIV.4.18 of~\cite{Tk3},
this algebra is a type II$_1$ factor,
whose unique tracial state $\ta_{\om}$
sends the image of a strongly $\om$-central sequence
$(a_l)_{l \in \N}$ in $N_{\om}$ to $\lim_{l \to \om} \ta (a_l).$
Let ${\overline{\af}}_{\om}$ be the induced automorphism of $N_{\om}.$
Theorem XIV.4.16 and Lemma XVII.2.2 of~\cite{Tk3} show that
${\overline{\af}}_{\om}^n$ is properly outer for all $n \neq 0.$
Apply Definition XVII.1.5 and Theorem XVII.1.6 of~\cite{Tk3}
to find \pj s
$f_0, f_1, \ldots, f_n \in N_{\om}$
such that $\sum_{k = 0}^n f_k = 1$
and
$\| {\overline{\af}}_{\om} (f_k) - f_{k + 1} \|_{2, \ta_{\om}} < \ep_0.$
By Theorem XIV.4.6(v) of~\cite{Tk3}, we can represent each
$f_k$ by a sequence $(f_{k, l})_{l \in \N}$ in $l^{\infty} (N)$
such that each $f_{k, l}$ is a \pj\  %
and $\sum_{k = 0}^n f_{k, l} = 1$ for every $l.$
For $0 \leq k \leq n$ and $a \in S$ we have
\[
\lim_{l \to \om} \ta ( [ a, \, f_{k, l}]^* [ a, \, f_{k, l}] )
    = \ta_{\om} ( [ a, \, f_{k}]^* [ a, \, f_{k}] )
    = 0,
\]
so there is a \nbhd\  $U$ of $\om$ in $\bt \N$ such that
$l \in U$ implies $\| [ a, \, f_{k, l}] \|_{2, \ta} < \ep_0$
for $a \in S$ and $0 \leq k \leq n.$
Similarly, for $l \in \N$ sufficiently close to $\om$ we have
$\| {\overline{\af}} (f_{k, l}) - f_{k + 1, \, l} \|_{2, \ta} < \ep_0$
for $0 \leq k \leq n - 1.$

Choose $l_0 \in \N$ for which the above estimates hold.
For $0 \leq k \leq n$ use Lemma~\ref{Z1} to find a \pj\  %
$g_k \in A$ such that
$\| g_k - f_{k, l_0} \|_{2, \ta}
  < \min \left( \frac{1}{2} \dt, \, \ep_0 \right).$
Then $\| g_j g_k \|_{2, \ta} < \dt$ for $j \neq k,$
so the choice of $\dt$ using Lemma~\ref{X2} provides
\mops\  $e_1, \ldots, e_n \in A$ such that
$\| e_k - g_k \|_{2, \ta} < \ep_0,$
and therefore $\| e_k - f_{k, l_0} \|_{2, \ta} < 2 \ep_0.$
Set $e_0 = 1 - \sum_{j = 1}^n e_j.$
Since $\sum_{j = 0}^{n} f_{j, l_0} = 1,$ we get
\[
\| e_0 - f_{0, l_0} \|_{2, \ta}
  \leq \sum_{j = 1}^n \| e_j - f_{j, l_0} \|_{2, \ta}
  < 2 n \ep_0.
\]
Now
\begin{align*}
\| \af (e_0) - e_1 \|_{2, \ta}
& \leq \| {\overline{\af}} (e_0 - f_{0, l_0}) \|_{2, \ta}
          + \| {\overline{\af}} (f_{0, l_0}) - f_{1, \, l_0} \|_{2, \ta}
          + \| f_{1, l_0} - e_1 \|_{2, \ta}  \\
& < 2 n \ep_0 + \ep_0 + 2 \ep_0
  < \ep.
\end{align*}
Similarly
$\| \af (e_j) - e_{j + 1} \|_{2, \ta} < 5 \ep_0 < \ep$
for $1 \leq j \leq n - 1.$
Finally, if $a \in S$ then
\[
\| [a, e_0] \|_{2, \ta}
  \leq 2 \| e_0 - f_{0, l_0} \|_{2, \ta} \cdot \| a \|
            + \| [a, f_{0, l_0}] \|_{2, \ta}
   < 2 \cdot 2 n \ep_0 + \ep_0
   < \ep,
\]
and similarly for $j \neq 0$ we have
$\| [a, e_j] \|_{2, \ta} < 2 \cdot 2 \ep_0 + \ep_0 < \ep.$
This completes the proof that outerness in the trace
representation implies the \tRp.

The converse follows from Proposition~2.3 of~\cite{Ks3}
and Corollary~4.6 of~\cite{OP1}.
\end{proof}

We can now give a version  of Kishimoto's result,
Theorem~2.1 of~\cite{Ks3}, giving conditions for
the Rokhlin property on a
simple unital AT~algebra with real rank zero and unique tracial state.

\begin{thm}\label{NTh}
Let $A$ be a simple separable unital C*-algebra with tracial
rank zero,
and suppose that $A$ has a unique tracial state $\ta.$
Let $\pi_{\ta} \colon A \to L (H_{\ta})$
be the Gelfand-Naimark-Segal representation associated with $\ta.$
Let $\af \in {\mathrm{Aut}}(A).$
Then the following condition are equivalent:
\begin{itemize}
\item[(1)]
$\af$ has the tracial Rokhlin property.
\item[(2)]
The automorphism of $\pi_{\ta} (A)''$ induced by $\af^n$
is outer for every $n \neq 0,$
that is, $\af^n$ is not weakly inner in $\pi_{\ta}$ for
any $n \neq 0.$
\item[(3)]
$C^*(\Z, A, \af)$ has a unique tracial state.
\item[(4)]
$C^*(\Z, A, \af)$ has real rank zero.
\end{itemize}
\end{thm}

\begin{proof}
The equivalence of~(1) and~(2) is Theorem~\ref{OuterImpTRP}.
The implication from~(1) to~(4) is Theorem~4.5 of~\cite{OP1}.
The implication from~(3) to~(2) is Proposition~2.3 of~\cite{Ks3}.

We prove that~(4) implies~(3).
Assume~(4).
Proposition~2.2 of~\cite{Ks3} implies that the restriction
map $T (C^* (\Z, A, \af)) \to T (A)$ is a bijection
to the set $T (A)^{\af}$ of $\af$-invariant tracial states of $A.$
Since $T(A)^{\af} = \{ \ta \},$
it follows that $C^*(\Z, A, \af)$ has a unique tracial state.
\end{proof}

\section{The tracial Rokhlin property for noncommutative
   Furstenberg transformations}\label{Sec:OutPf}

\indent
In this section, we prove that the automorphism $\af_{\te, \gm, d, f}$
of Definition~\ref{FDfn} has the \tRp\  when
$1, \, \te, \, \gm$ are linearly independent over $\Q.$
We do the same for the case $\te \in \R \setminus \Q,$ $\gm = 0,$
and $f$ constant.
Thus, the results of~\cite{OP1} apply to the crossed products by
these automorphisms.

A question that remains unanswered is whether these
automorphisms actually have the Rokhlin property.
An argument in Section~6 of~\cite{Ks1}
and an argument in~\cite{Ks2} proceed by first proving the
approximate Rokhlin property,
and then showing that, for the algebra and automorphism in question,
the approximate Rokhlin property implies the Rokhlin property.
Apparently, though, those arguments depend on approximate innerness.
Our automorphisms are not approximately inner because they
are nontrivial on $K_1.$

If $A$ is a simple \ca\  with unique tracial state $\ta,$
and $\af \in \Aut (A),$ we write ${\overline{\af}}$ for the
automorphism of $\pi_{\ta} (A)''$ determined by $\af.$
We also use the trace norm $\| \cdot \|_{2, \ta}$
described before Lemma~\ref{Q4}.

\begin{lem}\label{L2Sum}
Let $A$ be a separable unital \ca\  with a faithful tracial state $\ta.$
Let $(y_i)_{i \in I}$ be a family of unitaries such that
$\ta (y_i^* y_j) = 0$ for $i \neq j$
and whose linear span is dense in $A.$
Let $\pi \colon A \to L (H)$ be the
Gelfand-Naimark-Segal representation associated with $\ta.$
We identify $A$ with its image in $\pi (A)''.$
Then every $a \in \pi (A)''$ has a unique representation
as $a = \sum_{i \in I} \ld_i y_i,$
with convergence in $\| \cdot \|_{2, \ta}$ and whose coefficients
satisfy $\sum_{i \in I} | \ld_i |^2 = \| a \|_{2, \ta}^2.$
If $a$ is unitary then $\sum_{i \in I} | \ld_i |^2 = 1.$
\end{lem}

\begin{proof}
Let $\xi \in H$ be the standard cyclic vector for $\pi.$
Then one immediately checks that
$a \mapsto a \xi$ is an isometric linear map from
$\pi (A)''$ with $\| \cdot \|_{2, \ta}$ to $H,$
and that $(y_i \xi)_{i \in I}$ is an orthonormal basis for $H.$
The coefficients $\ld_i$ are determined by
$a \xi = \sum_{i \in I} \ld_i y_i \xi.$
\end{proof}

\begin{lem}\label{NCFIsOuterInTr}
Let $\te \in \R \setminus \Q$ and let $\gm \in \R.$
Let $\af \in \Aut (A_{\te})$ satisfy $\af (u) = \exp (2 \pi i \gm) u.$
If $\gm \not\in \Z + \te \Z$ then ${\overline{\af}}$ is outer.
\end{lem}

\begin{proof}
Let $\ta$ be the unique tracial state on $A_{\te},$
and identify $A_{\te}$ with its image in $\pi_{\ta} (A_{\te})''.$

Suppose ${\overline{\af}}$ is inner,
so ${\overline{\af}} = \Ad (w)$
for some unitary $w \in \pi_{\ta} (A_{\te})''.$
We apply Lemma~\ref{L2Sum} with $I = \Z \times \Z$
and $y_{m, n} = u^m v^n$ for $m, n \in \Z,$
and write
\[
w = \sum_{m, n \in \Z} \ld_{m, n} u^m v^n
\]
with convergence in $\| \cdot \|_{2, \ta}.$
Then, also with convergence in $\| \cdot \|_{2, \ta},$
\begin{align*}
\sum_{m, n \in \Z} e^{2 \pi i n \te} \ld_{m, n} u^{m + 1} v^n
 & = \sum_{m, n \in \Z} \ld_{m, n} u^m v^n u  \\
 & = w u
   = e^{2 \pi i \gm} u w
   = \sum_{m, n \in \Z} e^{2 \pi i n \gm} \ld_{m, n} u^{m + 1} v^n.
\end{align*}
We have $\ld_{m, n} \neq 0$ for some $m, \, n \in \Z,$
and uniqueness of the series representation then implies
$\exp (2 \pi i n \te) = \exp (2 \pi i \gm).$
\end{proof}

\begin{thm}\label{MainLemma}
Let $\te, \, \gm \in \R$ and suppose that
$1, \, \te, \, \gm$ are linearly independent over $\Q.$
Let $d \in \Z.$
Then the automorphism $\af = \af_{\te, \gm, d, f} \in \Aut (A_{\te}),$
of Definition~\ref{FDfn}, has the \tRp.
\end{thm}

\begin{proof}
It follows from Theorem~4 and Remark~6 of~\cite{EE}
and Proposition~2.6 of~\cite{LnTAF}
that $A_{\te}$ has tracial rank zero.
Also, $A_{\te}$ has a unique tracial state $\ta.$
Lemma~\ref{NCFIsOuterInTr} implies that ${\overline{\af}}^k$ is outer
on $\pi_{\ta} (A_{\te})''$ for all $k \neq 0.$
The theorem is therefore immediate from Theorem~\ref{OuterImpTRP}.
\end{proof}

Before stating the consequences, we recall a definition.

\begin{dfn}\label{OrdDetD}
Let $A$ be a unital \ca.
We say that the {\emph{order on \pj s over $A$ is determined by traces}}
if whenever $n \in \N$ and $p, q \in \Mi (A)$ are \pj s such that
$\ta (p) < \ta (q)$ for all $\ta \in T (A),$
then $p \precsim q.$
\end{dfn}

This is Blackadar's Second Fundamental Comparability Question
for $\Mi (A).$
See 1.3.1 in~\cite{Bl3}.

\begin{cor}\label{PropOfNCFT}
Let $\te, \, \gm \in \R$ and suppose that
$1, \, \te, \, \gm$ are linearly independent over $\Q.$
Let $d \in \Z.$
Let $\af_{\te, \gm, d, f} \in \Aut (A_{\te})$
be as in Definition~\ref{FDfn}.
Then:
\begin{itemize}
\item[(1)]
$C^* (\Z, \, A_{\te}, \, \af_{\te, \gm, d, f})$ is simple.
\item[(2)]
$C^* (\Z, \, A_{\te}, \, \af_{\te, \gm, d, f})$
has a unique tracial state.
\item[(3)]
$C^* (\Z, \, A_{\te}, \, \af_{\te, \gm, d, f})$ has real rank zero.
\item[(4)]
$C^* (\Z, \, A_{\te}, \, \af_{\te, \gm, d, f})$ has stable rank one.
\item[(5)]
The order on projections over
$C^* (\Z, \, A_{\te}, \, \af_{\te, \gm, d, f})$ is
determined by traces.
\item[(6)]
$C^* (\Z, \, A_{\te}, \, \af_{\te, \gm, d, f})$
satisfies the local approximation property of Popa~\cite{Pp}
(is a Popa algebra in the sense of Definition~1.2 of~\cite{Bn}).
\end{itemize}
\end{cor}

\begin{proof}
Part~(1)  follows from
Theorem~\ref{MainLemma} and Corollary~1.14 of~\cite{OP1}.
It is well known that $A_{\te}$ has real rank zero
(Theorem~1.5 of~\cite{BKR}; Remark~6 in Section~5 of~\cite{EE})
and stable rank one~(\cite{Pt1}),
and that the order on projections over $A_{\te}$ is
determined by traces (Corollary~2.5 of~\cite{Rf}).
(These also all follow from tracial rank zero;
see Theorem~3.4 of~\cite{LnTAF}
and Theorems~5.8 and~6.8 of~\cite{LnTTR}.)
Therefore Part~(2) follows from
Theorem~\ref{MainLemma} and Corollary~4.6 of~\cite{OP1},
Part~(3) follows from
Theorem~\ref{MainLemma} and Theorem~4.5 of~\cite{OP1},
Part~(4) follows from Theorem~\ref{MainLemma} and
Theorem~5.3 of~\cite{OP1},
Part~(5) follows from
Theorem~\ref{MainLemma} and Theorem~3.5 of~\cite{OP1},
and Part~(6) follows from
Theorem~\ref{MainLemma} and Corollary~4.7 of~\cite{OP1}.
\end{proof}

It is worth pointing out that one can verify the \tRp\  %
for the automorphisms above using methods related to those
of~\cite{Ph11},
without considering the automorphism in the trace representation
or using tracial rank zero.

We will show that the same methods give the \tRp\  for one other
kind of Furstenberg transformation on an irrational rotation algebra.

\begin{lem}\label{LinFIsOuterInTr}
Let $\te \in \R \setminus \Q.$
Let $\af \in \Aut (A_{\te})$ satisfy $\af (u) = u$
and $\af (v) = \ld u^d v$
with $\ld \in S^1$ and $d \in \Z \setminus \{ 0 \}.$
Then ${\overline{\af}}$ is outer.
\end{lem}

\begin{proof}
Let $\ta,$ $A_{\te} \subset \pi_{\ta} (A_{\te})'',$
${\overline{\af}} = \Ad (w),$
and $w = \sum_{m, n \in \Z} \ld_{m, n} u^m v^n$
be as in the proof of Lemma~\ref{NCFIsOuterInTr}.
Then, with convergence in $\| \cdot \|_{2, \ta},$
\begin{align*}
\sum_{m, n \in \Z} e^{2 \pi i n \te} \ld_{m, n} u^{m + 1} v^n
 & = \sum_{m, n \in \Z} \ld_{m, n} u^m v^n u  \\
 & = w u
   = u w
   = \sum_{m, n \in \Z} \ld_{m, n} u^{m + 1} v^n.
\end{align*}
Since $\te \not\in \Q,$
uniqueness of the series representation implies
$\ld_{m, n} = 0$ for $n \neq 0.$

Set $\bt_m = \ld_{m, 0}.$
Then, again with convergence in $\| \cdot \|_{2, \ta},$
\[
\sum_{m \in \Z} \bt_m u^m v
   = w v
   = \ld u^d v w
   = \sum_{m \in \Z} \ld \bt_m u^d v u^m
   = \sum_{m \in \Z} \ld e^{2 \pi i m \te} \bt_m u^{m + d} v.
\]
Uniqueness of the series representation implies
$| \bt_{m + d} | = | \bt_m |$ for all $m \in \Z.$
Since $\sum_{n = - \infty}^{\infty} | \bt_n |^2 = 1$ and $d \neq 0,$
this is a contradiction.
\end{proof}

In fact, if ${\overline{\af}}_{\te, 0, d, f}$ is inner,
then one can show that there are numbers $\bt_n \in \C$
for $n \in \Z$ with
$\sum_{n = - \infty}^{\infty} | \bt_n |^2 = 1$ such that
\[
e^{2 \pi i f (\zt)} \zt^d
  = \sum_{m, n \in \Z}
       \bt_m {\overline{\bt_n}} e^{- 2 \pi i n \te} \zt^{m - n}
\]
for all $\zt \in S^1.$

\begin{thm}\label{MainLemma2}
Let $\te \in \R \setminus \Q,$
let $\rh \in \R,$ and let $d \in \Z \setminus \{ 0 \}.$
Then the automorphism $\af = \af_{\te, 0, d, \rh} \in \Aut (A_{\te})$
of Definition~\ref{FDfn}
(with $\rh$ regarded as a constant function)
has the \tRp.
\end{thm}

\begin{proof}
The proof is the same as for Theorem~\ref{MainLemma}.
\end{proof}

We then get all the properties listed in Corollary~\ref{PropOfNCFT}
for the crossed products by these automorphisms as well.
By Lemma~\ref{ChangeOfVar}, the crossed product is also
the \ca\  of an ordinary smooth minimal Furstenberg transformation
on the torus, and in this case it is in fact known that the
\ca\  has tracial rank zero.
See~\cite{LP2}.

Theorem~\ref{MainLemma2} can also be proved by methods related to
those of~\cite{Ph11}, but with considerably more difficulty.

\begin{exa}\label{NonUErg}
We show that
Theorem~\ref{MainLemma} and Corollary~\ref{PropOfNCFT}(1) can fail
if $1, \, \te, \, \gm$ are not linearly independent over $\Q,$
even if the crossed product is simple,
$d = 1,$ and $\te$ and $\gm$ are both irrational.

By Theorem~2 (in Section~4) of~\cite{ILR} and the preceding discussion,
there is an irrational number
$\te \in \R$ (called $\af$ there) and a \cfn\  %
$g \colon S^1 \to \R$ (which in the notation of~\cite{ILR} is
$g (e^{2 \pi i t}) = f (t) - t$)
such that the function from $S^1$ to $S^1$
given by $\zt \mapsto e^{2 \pi i g (\zt)} \zt$ is a coboundary
with respect to the action of rotation by $e^{2 \pi i \te}$ on $S^1.$
Therefore the Furstenberg transformation
$h (\zt_1, \zt_2)
 = (\exp (2 \pi i \te) \zt_1, \, \exp (2 \pi i g (\zt_1)) \zt_1 \zt_2)$
is not ergodic with respect to Lebesgue measure.
(See the discussion before Proposition~1 in Section~4 of~\cite{ILR}.)
So $h$ is not uniquely ergodic.

Define $g_0 \colon S^1 \to \R$
by $g_0 (\zt) = g ( e^{- 2 \pi i \te} \zt) - \te.$
One checks that the automorphism $f \mapsto f \circ h^{-1}$
of $A_0 = C (S^1 \times S^1)$ is equal to
$\af_{0, \, - \te, \, - 1, \, - g_0}.$
It follows from Lemma~\ref{ChangeOfVar} that
$C^* (\Z, \, S^1 \times S^1, \, h)
 \cong C^* (\Z, \, A_{- \te}, \, \af_{- \te, \, 0, \, 1, \, g_0}).$
Accordingly,
the \ca\  $C^* (\Z, \, A_{- \te}, \, \af_{- \te, \, 0, \, 1, \, g_0})$
has more than one tracial state.
Since $A_{\te}$ has a unique tracial state,
it follows from Theorem~\ref{NTh} that
the automorphism $\af_{- \te, \, 0, \, 1, \, g_0}$
does not have the \tRp.

This version of the example does not have $\gm \not\in \Q,$
but by applying Lemma~\ref{Conj} with $k = 0$ and $l = - 1,$
we easily find a \cfn\  $k \colon S^1 \to \R$ such that
$C^* (\Z, \, A_{- \te}, \, \af_{- \te, \, 0, \, 1, \, g_0})
  \cong C^* (\Z, \, A_{- \te}, \, \af_{- \te, \, \te, \, 1, \, k}).$
So $C^* (\Z, \, A_{- \te}, \, \af_{- \te, \, \te, \, 1, \, k})$
does not have real rank zero, and
$\af_{- \te, \, \te, \, 1, \, k}$ does not have the \tRp.
\end{exa}

\section{C*-algebras of discrete subgroups of nilpotent Lie
  groups}\label{Sec:MW}

\indent
In this section, we consider the \ca s
$A_{\te}^{5, 3}$ and $A_{\te}^{5, 6}$ of~\cite{MW2},
which are the ``largest'' simple quotients of the \ca s of certain
discrete subgroups of five dimensional nilpotent Lie groups.
The corresponding results for $A_{\te}^{5, n}$ for $n = 1, 2, 4, 5$
follow from theorems already in the literature,
but those results don't apply
to $A_{\te}^{5, 3}$ and $A_{\te}^{5, 6}.$

\begin{prp}\label{MWAlg3}
For $\te \in \R \setminus \Q,$
the \ca\  $A_{\te}^{5, 3}$ of Section~3 of~\cite{MW2}
is isomorphic to a crossed product $C^* (\Z, B, \af)$
in which $B$ is the \ca\  %
of a smooth minimal Furstenberg transformation on $S^1 \times S^1$
and has a unique tracial state $\ta,$
and $\af$ has the \tRp.
\end{prp}

\begin{proof}
Let $\ld = \exp (2 \pi i \te).$
Then $A_{\te}^{5, 3}$ is the universal \ca\  generated by unitaries
$u, v, w, x$ satisfying the relations
\[
u v = x v u, \,\,\,\,\,\, u w = w u, \,\,\,\,\,\, u x = \ld x u,
\]
\[
v w = \ld w v, \,\,\,\,\,\, v x = x v, \andeqn w x = x w.
\]
(In~\cite{MW2}, the element $x$ is included in the relations
but not among the generators.
However, the first relation implies that $x$ is in the \ca\  %
generated by the other unitaries.)

We take $B$ to be the universal \ca\  generated by unitaries
$u, v, x$ satisfying the relations
\[
u v = x v u, \,\,\,\,\,\, u x = \ld x u, \andeqn v x = x v.
\]
We have
$B \cong C^* (\Z, \, S^1 \times S^1, \, h),$
where $h \colon S^1 \times S^1 \to S^1 \times S^1$ is
the smooth Furstenberg transformation given by
$h (\zt_1, \zt_2) = (\ld^{-1} \zt_1, \, \ld \zt_1^{-1} \zt_2),$
so that $h^{-1} (\zt_1 ,\zt_2) = (\ld \zt_1, \, \zt_1 \zt_2).$
Theorem~2.1 of~\cite{Fr} and the remark after it imply that
$h$ is uniquely ergodic and minimal.
Unique ergodicity implies that $B$ has a unique tracial state $\ta.$

One easily checks that there is a unique automorphism
$\af$ of $B$ satisfying
\[
\af (u) = u, \,\,\,\,\,\, \af (v) = \ld v,
\andeqn \af (x) = x,
\]
and that $C^* (\Z, B, \af) \cong A$ with the canonical unitary
of the crossed product being sent to $w^*.$

Let $\pi_{\ta}$ be the
Gelfand-Naimark-Segal representation associated to $\ta.$
We prove that for any $k \neq 0,$
the automorphism ${\overline{\af}}^k$ of $\pi_{\ta} (B)''$
induced by $\af^k$ is outer.
Suppose ${\overline{\af}}^k$ is inner,
so ${\overline{\af}}^k = \Ad (z)$
for some unitary $z \in \pi_{\ta} (B)''.$
We apply Lemma~\ref{L2Sum} with $I = \Z^3$
and $y_{l, m, n} = u^l v^m x^n$ for $l, m, n \in \Z,$
and write
\[
z = \sum_{l, m, n \in \Z} \ld_{l, m, n} u^l v^m x^n
\]
with convergence in $\| \cdot \|_{2, \ta}.$
Then, also with convergence in $\| \cdot \|_{2, \ta},$
\begin{align*}
\sum_{l, m, n \in \Z} \ld_{l, m, n} u^l v^m x^{n + 1}
 & = z x
   = x z
   = \sum_{l, m, n \in \Z} \ld_{l, m, n} x u^l v^m x^n   \\
 & = \sum_{l, m, n \in \Z} \ld_{l, m, n} \ld^{- l} u^l v^m x^{n + 1}.
\end{align*}
Uniqueness of the series representation and $\ld^{- l} \neq 1$
for $l \neq 0$
imply $\ld_{l, m, n} = 0$ for $l \neq 0.$
Therefore
\[
z = \sum_{m, n \in \Z} \ld_{0, m, n} v^m x^n.
\]
So $z$ commutes with $v,$ contradicting $\af^k (v) = \ld^k v.$
This shows that ${\overline{\af}}^k$ is outer.

It follows from~\cite{LP2} or~\cite{LhP}
that $A$ has tracial rank zero.
It now follows from Theorem~\ref{OuterImpTRP}
that $\af$ has the \tRp.
\end{proof}

\begin{cor}\label{PropOfMW3}
For $\te \in \R \setminus \Q,$
let $A_{\te}^{5, 3}$ be as in Section~3 of~\cite{MW2}.
Then:
\begin{itemize}
\item[(1)]
$A_{\te}^{5, 3}$ is simple.
\item[(2)]
$A_{\te}^{5, 3}$ has a unique tracial state.
\item[(3)]
$A_{\te}^{5, 3}$ has real rank zero.
\item[(4)]
$A_{\te}^{5, 3}$ has stable rank one.
\item[(5)]
The order on projections over
$A_{\te}^{5, 3}$ determined by traces.
\item[(6)]
$A_{\te}^{5, 3}$ satisfies
the local approximation property of Popa~\cite{Pp}
(is a Popa algebra in the sense of Definition~1.2 of~\cite{Bn}).
\end{itemize}
\end{cor}

\begin{proof}
Since the algebra $B$ in Proposition~\ref{MWAlg3} has tracial rank zero
(by~\cite{LP2}, or by \cite{LhP}),
it follows from Theorem~3.4 of~\cite{LnTAF} that $B$ has real rank zero,
and stable rank one,
and from Theorems~5.8 and~6.8 of~\cite{LnTTR} that
the order on projections over $B$ is determined by traces.
The proof is then the same as for Corollary~\ref{PropOfNCFT}.
\end{proof}

Parts~(1) and~(2) are proved in~\cite{MW2},
but the other properties are new.

\begin{rmk}\label{MWAlg3Gen}
The same results hold for the slightly more complicated versions
of these algebras which appear in~\cite{MW4}.
For $r_1, r_2, r_3, r_4, r_5 \in \Z$ let
$A_{\te}^{5, 3} (r_1, r_2, r_3, r_4, r_5)$ be the universal \ca\  %
generated by unitaries
$u, v, w, x$ satisfying the relations
\[
u v = \ld^{r_2} x^{r_1} v u, \,\,\,\,\,\,
 u w = \ld^{r_5} w u, \,\,\,\,\,\, u x = \ld^{r_3} x u,
\]
\[
v w = \ld^{r_4} w v, \,\,\,\,\,\, v x = x v, \andeqn w x = x w.
\]
(In~\cite{MW4}, the exponents $r_1, r_2, r_3, r_4, r_5$ are called
$\af, \bt, \gm, \dt, \ep.$)
Theorem~2 of~\cite{MW4} asserts that if
\[
r_1, r_3, r_4 > 0, \,\,\,\,\,\,
0 \leq r_5 \leq \tfrac{1}{2} {\mathrm{gcd}} (r_3, r_4),
\andeqn
0 \leq r_2 \leq \tfrac{1}{2} {\mathrm{gcd}} (r_1, r_3, r_4, r_3),
\]
then $A_{\te}^{5, 3} (r_1, r_2, r_3, r_4, r_5)$ is simple.
(These conditions are stated in Theorem~1 of~\cite{MW4}.
The last two conditions are really just normalization conditions,
and can be omitted in what follows.)

We apply the same analysis to these algebras.
The algebra $B$ of Proposition~\ref{MWAlg3}
is now the universal \ca\  generated by unitaries
$u, v, x$ satisfying the relations
\[
u v = \ld^{r_2} x^{r_1} v u,
 \,\,\,\,\,\, u x = \ld^{r_3} x u, \andeqn v x = x v.
\]
The \hme\  $h$ is now
\[
h (\zt_1, \zt_2)
 = (\ld^{-r_3} \zt_1, \, \ld^{r_1 r_3 - r_2} \zt_1^{- r_1} \zt_2).
\]
Since $r_1$ and $r_3$ are nonzero, it is still a
smooth uniquely ergodic minimal Furstenberg transformation.
The automorphism $\af$ is now determined by
\[
\af (u) = \ld^{- r_5} u, \,\,\,\,\,\, \af (v) = \ld^{r_4} v,
\andeqn \af (x) = x.
\]
Since $r_4 \neq 0,$
the proof that $\af$ has the \tRp\  is identical.
Thus, the conclusions of Corollary~\ref{PropOfMW3}
hold for these algebras as well.
Parts~(1) and~(2) are proved in~\cite{MW4}, but the others are new.
\end{rmk}

We now turn to $A_{\te}^{5, 6}.$

\begin{prp}\label{MWAlg6}
For $\te \in \R \setminus \Q,$
the \ca\  $A_{\te}^{5, 6}$ of Section~6 of~\cite{MW2}
is isomorphic to a crossed product $C^* (\Z, B, \af)$
in which $B$ is the \ca\  %
of a smooth minimal Furstenberg transformation on $S^1 \times S^1$
and has a unique tracial state $\ta,$
and $\af$ has the \tRp.
\end{prp}

\begin{proof}
Let $\ld = \exp (2 \pi i \te).$
Then $A_{\te}^{5, 6}$ is the universal \ca\  generated by unitaries
$u, v, w, x$ satisfying the relations
\[
u v = w v u, \,\,\,\,\,\, u w = x w u, \,\,\,\,\,\, u x = \ld x u,
\]
\[
v w = \ld w v, \,\,\,\,\,\, v x = x v, \andeqn w x = x w.
\]
(In~\cite{MW2}, the elements $w$ and $x$ are included in the relations
but not among the generators.
However,
the first two relations imply that $w$ and $x$ are in the \ca\  %
generated by $u$ and $v.$)

We take $B$ to be the universal \ca\  generated by unitaries
$u, w, x$ satisfying the relations
\[
u w = x w u, \,\,\,\,\,\, u x = \ld x u, \andeqn w x = x w.
\]
As in the proof of Proposition~\ref{MWAlg3},
this \ca\  is the crossed product by a uniquely ergodic
minimal smooth Furstenberg transformation on $S^1 \times S^1.$
One easily checks that there is a unique automorphism
$\af$ of $B$ satisfying
\[
\af (u) = w^* u, \,\,\,\,\,\, \af (w) = \ld w,
\andeqn \af (x) = x,
\]
and that $C^* (\Z, B, \af) \cong A$ with the canonical unitary
of the crossed product being sent to $v.$

Let $\ta$ be the unique tracial state on $B,$
and let $\pi_{\ta}$ be the associated
Gelfand-Naimark-Segal representation.
We prove that for $k \neq 0,$
the automorphism ${\overline{\af}}^k$ of $\pi_{\ta} (B)''$
induced by $\af^k$ is outer.
Suppose ${\overline{\af}}^k$ is inner,
so ${\overline{\af}}^k = \Ad (z)$
for some unitary $z \in \pi_{\ta} (B)''.$
We apply Lemma~\ref{L2Sum} with $I = \Z^3$
and $y_{l, m, n} = u^l w^m x^n$ for $l, m, n \in \Z,$
and write
\[
z = \sum_{l, m, n \in \Z} \ld_{l, m, n} u^l w^m x^n
\]
with convergence in $\| \cdot \|_{2, \ta}.$
Then, also with convergence in $\| \cdot \|_{2, \ta},$
\begin{align*}
\sum_{l, m, n \in \Z} \ld_{l, m, n} u^l w^m x^{n + 1}
 & = z x
   = x z
   = \sum_{l, m, n \in \Z} \ld_{l, m, n} x u^l w^m x^n   \\
 & = \sum_{l, m, n \in \Z} \ld_{l, m, n} \ld^{- l} u^l w^m x^{n + 1}.
\end{align*}
Uniqueness of the series representation and $\ld^{- l} \neq 1$
for $l \neq 0$
imply $\ld_{l, m, n} = 0$ for $l \neq 0.$
Therefore
\[
z = \sum_{m, n \in \Z} \ld_{0, m, n} w^m x^n.
\]
So $z$ commutes with $w,$ contradicting $\af^k (w) = \ld^k w.$
This shows that ${\overline{\af}}^k$ is outer.

It follows from ~\cite{LP2} or~\cite{LhP}
that $A$ has tracial rank zero.
It now follows from Theorem~\ref{OuterImpTRP}
that $\af$ has the \tRp.
\end{proof}

\begin{cor}\label{PropOfMW6}
For $\te \in \R \setminus \Q,$
let $A_{\te}^{5, 6}$ be as in Section~6 of~\cite{MW2}.
Then:
\begin{itemize}
\item[(1)]
$A_{\te}^{5, 6}$ is simple.
\item[(2)]
$A_{\te}^{5, 6}$ has a unique tracial state.
\item[(3)]
$A_{\te}^{5, 6}$ has real rank zero.
\item[(4)]
$A_{\te}^{5, 6}$ has stable rank one.
\item[(5)]
The order on projections over
$A_{\te}^{5, 6}$ determined by traces.
\item[(6)]
$A_{\te}^{5, 6}$ satisfies the local approximation property of Popa.
\end{itemize}
\end{cor}

\begin{proof}
The proof is the same as for Corollary~\ref{PropOfMW3}.
\end{proof}

As for $A_{\te}^{5, 3},$
Parts~(1) and~(2) are proved in~\cite{MW2},
but the other properties are new.

\end{document}